\documentclass[a4paper, 11pt]{article}

\usepackage{amssymb,amsthm,amsmath,amsfonts,amsbsy}
\usepackage[affil-it, auth-lg]{authblk}
\usepackage{scrextend}
\usepackage[T1]{fontenc}
\usepackage{anyfontsize}

\usepackage{mathrsfs}
\usepackage{upgreek}
\usepackage{pifont, dsfont}
\usepackage{srcltx}
\usepackage[usenames,dvipsnames]{color}
\usepackage{url}
\usepackage[colorlinks,citecolor=dblue,urlcolor=lblue,
			linkcolor=dblue,breaklinks]{hyperref}
\usepackage[shortlabels]{enumitem}
\usepackage{titlesec}
\usepackage{calc}
\usepackage{subfig}
\usepackage{multicol}
\usepackage{rotating}
\usepackage{overpic}
\usepackage[round]{natbib}
\usepackage{caption}
\usepackage{fancyhdr}

\usepackage[width=6.25in,height=9.25in]{geometry}

\newcommand{\doi}[1]{\scshape{doi}:\,\tt\href{http://dx.doi.org/#1}
{\nolinkurl{#1}}\normalfont}

\titleformat*{\section}{\large\bf}
\titleformat*{\subsection}{\large\sl}

\theoremstyle{plain}

\newtheorem{thm}{\bf T\footnotesize HEOREM\normalsize}[section]

\newtheorem{lem}[thm]{\bf L\footnotesize EMMA\normalsize}
\newtheorem{prop}[thm]{\bf P\footnotesize ROPOSITION\normalsize}

\theoremstyle{definition}
\newtheorem{defn}[thm]{\bf D\footnotesize EFINITION\normalsize}
\newtheorem{rmk}[thm]{\bf R\footnotesize EMARK\normalsize}
\newtheorem{example}[thm]{\bf E\footnotesize XAMPLE\normalsize}

\numberwithin{equation}{section}

\newcommand{\thmref}[1]{\hyperref[#1]{Theorem~\ref*{#1}}}
\newcommand{\lemref}[1]{\hyperref[#1]{Lemma~\ref*{#1}}}
\newcommand{\defnref}[1]{\hyperref[#1]{Definition~\ref*{#1}}}
\newcommand{\egref}[1]{\hyperref[#1]{Example~\ref*{#1}}}
\newcommand{\chpref}[1]{\hyperref[#1]{Chapter~\ref*{#1}}}
\newcommand{\appref}[1]{\hyperref[#1]{Appendix~\ref*{#1}}}
\newcommand{\rmkref}[1]{\hyperref[#1]{Remark~\ref*{#1}}}
\newcommand{\secref}[1]{\hyperref[#1]{Section~\ref*{#1}}}
\newcommand{\assref}[1]{\hyperref[#1]{Assumption~\ref*{#1}}}
\newcommand{\propref}[1]{\hyperref[#1]{Proposition~\ref*{#1}}}
\newcommand{\cororef}[1]{\hyperref[#1]{Corollary~\ref*{#1}}}
\newcommand{\figref}[1]{\hyperref[#1]{Figure~\ref*{#1}}}
\newcommand{\itmref}[1]{\hyperref[#1]{(\ref*{#1}}}
\renewcommand{\eqref}[1]{\hyperref[#1]{(\ref*{#1})}}
\newcommand{\ppref}[1]{\hyperref[#1]{p.\pageref*{#1}}}

\newcommand{\Cf}{\mathcal{C}}
\newcommand{\cl}[1]{\mathrm{cl}(#1)}
\newcommand{\dif}{\mathrm{d}}
\newcommand{\E}{\mathbb{E}}
\newcommand{\eps}{\varepsilon}
\newcommand{\half}{\dfrac{1}{2}}
\newcommand{\ind}{\mathds{1}}
\newcommand{\pb}{\mathbb{P}}
\newcommand{\pd}[2]{\dfrac{\partial #1}{\partial #2}}
\newcommand{\pdd}[2]{\dfrac{\partial^{2}#1}{\partial #2^{2}}}
\newcommand{\poten}[1]{\mathrm{U}^{#1}}
\newcommand{\R}{\mathbb{R}}
\newcommand{\sgn}{\mathrm{sgn}}
\newcommand{\wti}{\widetilde}
\newcommand{\mmu}{\boldsymbol{\upmu}}
\newcommand{\mtau}{\uptau}
\newcommand{\mtheta}{\uptheta}
\newcommand{\mrho}{\uprho}

\newcommand{\SEP}{\ref{eq:SEP}}
\newcommand{\OBS}{\ref{eq:OBS}}
\newcommand{\mSEP}{\ref{eq:mSEP}}
\newcommand{\mOBS}{\ref{eq:mOBS}}

\DeclareMathAlphabet{\mathcal}{OMS}{cmsy}{m}{n}

\definecolor{dblue}{RGB}{0, 72,128}
\definecolor{lblue}{RGB}{0,128,192}

\allowdisplaybreaks

\makeatletter
\renewcommand*{\@fnsymbol}[1]{\ensuremath{\ifcase#1\or 
	\text{\,\footnotesize\ding{73}}\or \text{\,\small$\ast$}\or \ddagger\or
    \mathsection\or \mathparagraph\or \|\or **\or \dagger\dagger
    \or \ddagger\ddagger \else\@ctrerr\fi
    }}
\makeatother

\title{\vspace{-.5em}Minimal Root's Embeddings for General Starting
		and Target Distributions\thanks{\,This research is financially 
			supported by the National Natural Science Foundation of China 
			(No.\,11601306) and the Research Foundation of Jinan 
			University (No.\,21617415 \& 88018052). This paper was 
			partially completed while the author visited Sapienza University 
			of Rome as a post-doctoral fellow (SECS-S/06, 2013--2015). 
			\emph{This is an electronic reprint of the manuscript accepted
			by \href{https://www.elsevier.com/locate/spa}{\emph{Stochastic 
			Processes and their Applications}}, \doi{10.1016/j.spa.2019.01.009}. 
			\it This reprint differs from the published paper in pagination and 
			typographic detail.}\vspace*{0.3em}}
		}

\author[ ]{Jiajie Wang\thanks{\,\footnotesize 
						Correspondence to: SZTC, Jinan University, 
						Shenzhen, Guangdong, 518053, P.R. China.\newline
						\hspace*{2em}E-mail addresses:\ \ \upshape
						\href{mailto:wang_jj@sz.jnu.edu.cn}
							{\nolinkurl{wang_jj@sz.jnu.edu.cn}},\ \ 			
						\href{mailto:jiajie.wang@bath.edu}
							{\nolinkurl{jiajie.wang@bath.edu}}.
						}}

\affil[ ]{\small\hspace{0.4pt}SZTC, Jinan University, Shenzhen, Guangdong, 
		518053 - P.R. China.\vspace{2pt}}
\affil[ ]{\small\hspace{0.4pt}MEMOTEF, Sapienza University of Rome, Rome, 
		00161 - Italy.\vspace*{-9pt}}

\date{}

\begin{document}

\linespread{1.1}

\maketitle

\fancyhead{}
\fancyfoot{}
\fancyhead[L]{\scriptsize
			\copyright\ 2019. This manuscript version is made available
			under the CC-BY-NC-ND 4.0 license.\newline
			\hspace*{0.1em}\url{http://creativecommons.org/licenses/by-nc-nd/4.0/}
			}
\renewcommand{\headrulewidth}{0.00pt}

\vspace*{-15pt}
\thispagestyle{fancy}


	
\begin{abstract}

	\noindent 
	Most results regarding Skorokhod embedding problems (\SEP) so far rely 
	on the assumption that the corresponding stopped process is uniformly
	integrable, which is equivalent to the convex ordering condition
	$\poten{\mu}\leq \poten{\nu}$ when the underlying process is a local
	martingale. In this paper, we study the existence, construction of
	Root's solutions to \SEP, in the absence of this convex ordering
	condition. We replace the uniform integrability condition by the 
	minimality condition (\citealp{Monroe:1972a}), as the criterion of
	``good'' solutions. A sufficient and necessary condition (in terms of
	local time) for minimality is given. We also discuss the optimality of
	such minimal solutions. These results extend the generality of the
	results given by \citet{CoxWang:2013a} and 
	\citet{GassiatOberhauserReis:2015}. At last, we extend all the results
	above to multi-marginal embedding problems based on the work of 
	\citet{CoxOblojTouzi:2018}. 
	
	\vspace*{6pt}
	
\end{abstract}


\vspace*{10pt}

\normalsize

\linespread{1.2}
	

\section{Introduction}
\label{sec:introduction}
	
Given a stochastic process $X$ on a filtered probability space $(\Omega,
\mathcal{F}, \mathscr{F}=\{\mathcal{F}_t\}_{t\geq0},\pb)$, and a 
distribution $\mu$ on the state space of $X$, the Skorokhod embedding 
problem is to find a stopping time $\tau$ such that $X_\tau\sim\mu$. This 
problem was initially proposed by \citet{Skorokhod:1965}. 
	
Under the classical setting where $X$ is a Brownian motion starting at $0$ 
and the target distribution $\mu$ has zero mean and finite variance, there 
is a rich literature regarding this problem, for example, 
\citet{Dubins:1968}, \citet{Root:1969}, \citet{Rost:1971},
\citet{ChaconWalsh:1976}, \citet{AzemaYor:1979a}, \citet{Vallois:1983}, 
\citet{ChaconRene:Thesis}, \citet{Perkins:1985}, etc. We will not state 
them one by one in details. Instead, we refer curious readers to the 
survey paper \citet{Obloj:2004}.
	
Most of the results above can be generalized to the cases where the 
underlying process is a diffusion process with a general starting 
distribution. In this paper, we denote such embedding problem by 
\SEP$(\sigma,\nu,\mu)$:
\begin{equation}\tag{$\mathrm{SEP}$}\label{eq:SEP}
	\begin{aligned}
		\textit{Given $X_0\sim\nu$, \,to
			find a stopping time $\tau$ such that 
			$X_\tau\sim\mu$\,,}\\[4pt]
		\textit{where $X$ satisfies 
			$\dif X_t=\sigma(X_t)\dif W_t$\,.}
	\end{aligned}
\end{equation}
	
However, the results mentioned above are concerned with the cases where 
the embeddings are namely \emph{UI stopping times}. Here, a stopping time 
$\tau$ is a UI stopping time if the corresponding stopped process $X^{\tau}
:=\{X_{t\wedge\tau}\}_{t\geq0}$ is uniformly integrable, otherwise we call 
$\tau$ a \emph{non-UI stopping time}.
	
When the underlying process is a continuous local martingale 
\citet[Prop.\,8.1]{Obloj:2004} shows that there exists a UI embedding for 
\SEP$(\sigma,\nu,\mu)$ if and only if \emph{the convex ordering condition} 
holds
\begin{equation}\label{eq:regular}
	\poten{\nu}(x)\ \geq\ \poten{\mu}(x)\ >\ -\infty\,,
	\quad\textit{for all}\ \ x\in\R,
\end{equation}
where the function $\poten{\mu}$ is called \emph{the potential of $\mu$} 
(\citealp{Chacon:1977}):
\begin{align*}
	\poten{\mu}(x)\ :=\ -\E^{Y\sim \mu}\left[\,\big|Y-x\big|\,\right]
	\ =\ -\int_\R|y-x|\,\mu(\dif y)\,.
\end{align*}
We say that \emph{$\nu\preceq\mu$ in convex order} if \eqref{eq:regular}
holds.
	
\emph{In this paper, we are concerned with \ref{eq:SEP}\,$(\sigma,\nu,\mu)$ 
in the absence of convex ordering condition \eqref{eq:regular}}. In such 
circumstances we cannot expect the corresponding embedding to be a UI 
stopping time.

For example,  suppose that the initial distribution is the Dirac measure
$\nu=\delta_0$ and the target is $\mu=\delta_1$. The mean values of $\nu$ 
and $\mu$ do not agree, and then \eqref{eq:regular} fails. The hitting time 
$H_1=\inf\{t\geq0:\,W_t=1\}$ is an embedding for \SEP$(\nu,\mu)$ but 
obviously it is \emph{not} a UI stopping time. Another example is that 
$\nu=(\delta_1+\delta_{-1})/2$ and $\mu=\delta_0$. The mean values agree, 
but \eqref{eq:regular} fails as $\poten{\nu}\leq \poten{\mu}$. The hitting 
time $H_0=\inf\{t\geq0: \,W_t=0\}$ is a non-UI embedding for 
\SEP$(\nu,\mu)$.

As presented above, in the absence of \eqref{eq:regular}, we cannot restrict 
our attention to UI stopping times for embeddings. Instead, we may pose some 
other restrictions.  For example, \citet{PedersenPeskir:2001} pose an 
integrability condition on the maximum of the scale function of $X$ as the 
replacement of UI condition. After that, \citet{CoxHobson:2006} propose 
another criterion on stopping times, which was initially introduced by  
\citet{Monroe:1972a}:
\begin{defn}[Minimal stopping time]\label{defn:minimal}
	A stopping time $\tau$ for the process $X$ is minimal if whenever 
	$\theta\leq \tau$ is a stopping time such that $X_\theta$ and $X_\tau$ 
	have the same distribution then $\tau=\theta$, a.s..
\end{defn}
	
According to the definition, minimal stopping times could be a natural 
choice for ``good'' solutions of the embedding problem in a general context. 
For example, as stated in \citet[Sect.\,4.2]{Hobson:2011}, there exists a 
trivial solution for \SEP\ in the general cases\ ---\ simply run the process 
$X$ until it firstly hits the mean of $\mu$, and thereafter can use any 
regular embedding mentioned above. The embeddings constructed in this way 
are always minimal stopping times, see \citet{CoxHobson:2006}.

\citet{CoxHobson:2006} have made significant effort in the study of minimal 
stopping times for the Brownian motion starting at $0$. A group of necessary 
and sufficient conditions for minimality is given. After that, 
\citet{Cox:2008} extends the previous results to the cases of general 
starting distributions.  Thanks to these results, some well-known embeddings 
have been extended to the cases in which \eqref{eq:regular} fails, such as 
Chacon-Walsh's embedding, Az\'ema-Yor's embedding, Vallois' embedding.

In this work we are concerned with embeddings of Root's type which was
initially proposed by \citet{Root:1969}. Formally, suppose that $W$ is a 
Brownian motion starting at zero and the target distribution is a centred 
distribution with finite second moment, the Skorokhod embedding problem 
admits a solution which is the first hitting time of the joint process 
$(t,W_t)$ of a called \emph{Root's barrier}:
\begin{defn}[Root's barrier] \label{defn:barrier} 
	A closed subset $B$ of $[0,+\infty]\times[-\infty,+\infty]$ 
	is a \emph{Root's barrier} if
		
	\vspace{3pt}\noindent 
	a). \ $(+\infty,x)\in B$ \ if \ $x\in[-\infty,+\infty]$; \hfill\hfill
	b). \ $(t,\pm\infty)\in B$ \ if \ $t\in[0,\infty]$;\hfill\,
		
	\vspace{3pt}\noindent
	c). \ if $(t,x)\in B$, then $(s,x)\in B$ \ whenever $s>t$.
\end{defn}	
	
There have been a number of important contributions concerning Root's 
barriers (given that \eqref{eq:regular} holds). An immediately subsequent 
paper \citet{Loynes:1970} shows some elementary analytical properties of 
Root's barriers. Further, by posing the definition \emph{regular barrier}, 
the uniqueness of Root's embedding is given in this paper.

Another important paper regarding Root's construction is \citet{Rost:1976} 
which vastly extends the generality of Root's existence result. More 
importantly, R\"ost firstly proved the optimality of Root's embedding, which 
was conjectured by \citet{Kiefer:1972}, in the sense of \emph{minimal 
residual expectation} (\ref{eq:mre}, for short):
\begin{equation}\label{eq:mre}\tag{m.r.e.}
	\begin{aligned}
		\textit{Amongst all solutions of \,\SEP$\,(\sigma,\nu,\mu)$, 
			\,the Root's solution}&\\[4pt]
		\textit{minimises $\E^\nu\big[(\tau-t)^+\big]$ simultaneously 
			for all $t>0$}&.
	\end{aligned}
\end{equation}

\citet{Dupire:2005} proposes the connections among Root's embeddings, PDE 
and robust pricing problem for variance options. Enlightened by his idea, we 
derive the construction  of Root's embeddings using variational inequalities 
(given that \eqref{eq:regular} holds) in \citet{CoxWang:2013a}. We also 
propose the conjecture that, by slightly changing the terminal condition in 
our variational inequalities, this construction method could be extended to 
the cases where \eqref{eq:regular} fails 
(\citealp[Rmk.\,4.5]{CoxWang:2013a}). In the same paper, an alternative 
proof of \ref{eq:mre} property is given, which has an important application 
for the construction of sub-hedging strategies in the financial context. 
Later, using PDE techniques, \citet{GassiatOberhauserReis:2015} describe 
Root's embedding in terms of viscosity solutions of obstacle problems, and 
give a rigorous proof of the existence of Root's embedding given
\eqref{eq:regular}; using method from optimal transport,
\citet{BeiglbockCoxHuesmann:2017} show same existence and optimality results 
of Root's barriers. A more recent paper, \citet{CoxOblojTouzi:2018}, 
discusses the multi-marginal \mSEP, which is to find an increasing sequence 
of stopping times embedding the given multiple target distributions (in 
convex order) in sequence. They construct the UI solution of Root's type to 
the multi-marginal \mSEP\ via iterated optimal stopping problems. The 
optimality of such solutions is also given in their work.

In this work, we will extend the generality of the construction given by
\citet{CoxWang:2013a} and \citet{GassiatOberhauserReis:2015} to the cases 
without convex ordering condition \eqref{eq:regular}. On the other hand, 
thanks to the rich results given in \citet{CoxHobson:2006} and 
\citet{Cox:2008}, it will turn out that we can characterize minimal stopping 
times by the local times of the corresponding stopped process ($\E^\nu
\big[ L^x_\tau\big]$). This characterization then ensure that we can 
construct a minimal Root's embedding via an obstacle problem with proper 
boundary condition. Using the result about minimality, we then can discuss 
optimality of minimal Root's solutions (among all minimal solutions). After 
that, based on the work of \citet{CoxOblojTouzi:2018}, it turns out that one 
can construct Root's solution to multi-marginal \mSEP\ via iterated obstacle 
problems even when convex ordering condition fails. Moreover, we define the 
minimality for a sequence of stopping times, and tell when the solution to a 
multi-marginal \mSEP\ is ``minimal''. 

The paper will therefore proceed as follows: in \secref{sec:preliminaries}, 
we review some early results about Root's barriers. In 
\secref{sec:construction}, the existence result and the construction of 
Root's barrier for general starting and target distributions are given. In 
\secref{sec:minimality}, we study the potentials of the corresponding 
stopped process (and their limit), and obtain a necessary and sufficient 
condition for a Root's stopping time to be minimal. In 
\secref{sec:optimality}, we consider the optimality of non-UI Root's 
embeddings in the sense of \emph{maximal principal expectation}, which can 
be regarded as the generalization of minimal residual expectation 
(\ref{eq:mre}). In \secref{sec:multi_margin}, we extend all the results 
(construction, minimality, optimality) to the embedding problems with 
multi-marginal distributions.


\section{Preliminaries: Root's barriers for regular cases}
\label{sec:preliminaries}
	
We firstly review the previous results regarding Root's embeddings, which 
are useful throughout this work.

It was shown in \citet[Prop.\,3]{Loynes:1970} that the set $B$ defined in 
\defnref{defn:barrier} can be represented as a closed set bounded below by 
a lower semi-continuous function $R:\R\rightarrow[0,+\infty]$, i.e. 
$B=\big\{(t,x):\ t\geq R(x)\big\}$. This representation has been helpful in 
the characterization of the law of the stopped process $X^\tau$. 
Additionally, in the rest of this paper, we will say that a barrier is 
either a closed set described in \defnref{defn:barrier}, or equivalently its 
complement:
\begin{align*}
	D\ =\ \big\{(t,x):\ 0<t<R(x)\big\}
		\ =\ \left(\R_+\times\R\right)\backslash B.
\end{align*}
The corresponding stopping time is denoted by
\begin{equation*}
	\tau_D\ :=\ \inf\big\{\,t>0:\,(t,X_t)\notin D\,\big\}\ =\ 
		\inf\big\{\,t>0:\,t\geq R(X_t)\,\big\}.
\end{equation*}
Moreover, \citet[Prop.\,1]{Loynes:1970} says that, 
\begin{align*}
	\text{ for a Root's stopping $\tau_D$, \ either \ 
		$\pb[\tau_D<\infty]=1$ \ or \ $\pb[\tau_D=\infty]=1$.}
\end{align*}
	
As a straightforward result of this proposition, when $X_{\tau_D}\sim\mu$ 
where $\mu$ is integrable, $\tau_D$ is finite almost surely. 
	
The following properties are given in \citet{CoxWang:2013a}, which enable 
us to characterize the behaviour of the path of corresponding stopped 
process.
\begin{prop}\label{prop:path}
	Suppose $X$ is a continuous process. Given a Root's barrier $D$ and 
	the corresponding stopping time is denoted by $\tau_D$, then
	\begin{enumerate}[i), leftmargin = 2.1em, labelsep = .9em ]
		\item\label{itm:insideD}
			if $(t,x)\in D$, $\ \pb^\nu\big[X_{t\wedge\tau_D}\in\dif x\big]
				\,=\,\pb^\nu\big[X_{t}\in\dif x,\,t<\tau_D\big]$;
		\item\label{itm:outsideD}
			if $(t,x)\notin D$, $\ \pb^\nu\left[L^x_{t\wedge\tau_D}
				\,=\,L^x_{\tau_D}\right]\,=\,1$.
	\end{enumerate}
\end{prop}
	
These properties are local properties and do not rely on the integrability 
of the stopped process $X^\tau$, so they remain true even when $X^\tau$ is 
not uniformly integrable.

Denote the potential of the stopped process by $u(t,x)\,:=\,-\E^\nu
\big|x-X_{t\wedge\tau_D}\big|$, then, according to \citet{CoxWang:2013a}
(see also \citet{GassiatOberhauserReis:2015}), $u$ is of the class
$\Cf^0(\mathbb{R}_+\times\mathbb{R})\cap\Cf^{2,1}(D)$, and satisfies
\begin{subequations}
	\begin{align}\label{eq:general}
		\mathrm{L}u\,:=\,\pd{u}{t}-\dfrac{\sigma^2}{2}\pdd{u}{x}\,=\,0
			\quad\text{on}~~D;
			\qquad u(0,\cdot)\,=\,\poten{\nu}(\cdot)\quad\text{on}~~\R.
	\end{align}
	Moreover, if $\tau_D$ is a UI stopping time such that $X_{\tau_D}\sim
	\mu$, then 
	\begin{align}\label{eq:onlyUI}
		u(t,x)\,=\,\poten{\mu}(x),~~\text{if}~~(t,x)\notin D;
			\quad u(t,x)\,\longrightarrow\,\poten{\mu}(x),
			~~\text{as}~~t\rightarrow\infty.
	\end{align}
\end{subequations}
Note that the UI condition implies that $\poten{\nu}\geq \poten{\mu}$ 
everywhere on $\R$.	
	
In \citet{CoxWang:2013a}, we consider Root's embeddings for homogeneous 
diffusions, i.e. $\sigma(t,x)\equiv\sigma(x)$. Suppose that 
\eqref{eq:regular} holds, using \eqref{eq:general}-\eqref{eq:onlyUI},	
we construct a one-to-one correspondence between Root's stopping times and 
strong solutions to variational inequalities. Later, using the notion of 
viscosity solutions, \citet{GassiatOberhauserReis:2015} extend the result 
to more general cases.
\begin{thm}[\citealp{GassiatOberhauserReis:2015}]\label{thm:Gassiat}
	Assume that $(\nu,\mu)$ satisfies \eqref{eq:regular}, and $\sigma$ 
	satisfies that the following regular conditions:
	\begin{align*}
		&\begin{aligned}
			\text{there~exists}~~L>0,~~\text{s.t.}~~
				\forall\,t\geq 0,~~x,y\in\R,&\\[4pt]
			|\sigma(t,x)-\sigma(t,y)|<L|x&-y|,~~~|\sigma(t,x)|<L(1+|x|)\,;
		\end{aligned}\\[8pt]
		&\begin{aligned}
			\text{for~each~compact}~~K\subset\big\{x:\poten{\nu}(x)>
				\poten{\mu}(x)\big\},&\\[4pt]
			\exists~C_K>0,~~\text{s.t.}~~\forall\,t\geq 0,~~&x\in K,
				~~\sigma(t,x)\,\geq\,C_K\,>\,0\,.
		\end{aligned}
	\end{align*}
	Further, let $\tau_D$ be a UI Root's solution to 
	\SEP$(\sigma,\nu,\mu)$, and the function $u(t,x)$ be a viscosity 
	solution to the following obstacle problem
	\begin{equation*}
		\min\big\{\,\mathrm{L}u,\ u-\poten{\mu}\,\big\}\,=\,0,
			\quad u(0,\cdot)\,=\,\poten{\nu}(\cdot)\,,
			\quad \lim_{t\to\infty}u(t,\cdot)\,=\,\poten{\mu}(\cdot)
	\end{equation*}
	Then $\,u(t,x)=-\E^\nu\big|x-X_{t\wedge\tau_D}\big|$ \ and \
	$D=\big\{(t,x):\, u(t,x)>\poten{\mu}(x)\big\}$.
\end{thm}	
	
As stated in \secref{sec:introduction},$\tau_D$ is non-UI when 
\eqref{eq:regular} fails, and then \eqref{eq:onlyUI} does not hold any 
longer. Consequently, the results of \cite{CoxWang:2013a} and
\citet{GassiatOberhauserReis:2015} are not available. However, since 
\eqref{eq:general} still holds, in order to construct Root's embedding, 
we only need to find a more general version of \eqref{eq:onlyUI} --- it is 
the starting point of this work.


\section{Existence and construction of Root's embeddings}
\label{sec:construction}
	
Given $\poten{\mu}\leq\poten{\nu}$, \SEP$(\sigma,\nu,\mu)$ admits a UI 
Root's solution, and we can construct this solution via an obstacle problem
(\thmref{thm:Gassiat}). However, when \eqref{eq:regular} fails, we cannot 
even be sure if the Root's embedding exists. From now on, we are concerned 
with the existence and construction of Root's embedding in such general 
cases.	
	
First of all, let $\nu$ and $\mu$ be two probability distributions on $\R$, 
and define
\begin{align}\label{eq:OBS_setting}
	\begin{aligned}
		u_0(x)\,=\,\poten{\nu}(x),\quad\bar{u}(x)\,=\,\poten{\mu}(x)-C,
			\quad\text{for}~x\in\R,\\[6pt]
		\text{where $C>0$ is a constant s.t.~$u_0\geq\bar{u}$~everywhere}.
	\end{aligned}
\end{align}
We assume that the diffusion coefficient $\sigma$ satisfies the regular 
conditions:
\begin{align}
	&\label{eq:Lip_LG}
	\begin{aligned}
		\text{there exists}~L>0,~~\text{s.t.}~~\forall\,x,y\in\R,&\\[4pt]
		|\sigma(x)-\sigma(y)|<&\,L|x-y|,~~|\sigma(x)|<L(1+|x|)\,;
	\end{aligned}\\[8pt]
	&\label{eq:ellipticity}
	\begin{aligned}
		\text{for each compact}~K\subset&\,
			\big\{x:u_0(x)>\bar{u}(x)\big\},\\[4pt]
		\exists&~C_K>0,~~\text{s.t.}~~\forall\,x\,\in K,
			~~\sigma(x)\,\geq\,C_K\,>\,0\,.
	\end{aligned}
\end{align}
Consider the obstacle problem \OBS$(\sigma,u_0,\bar{u})$:
\begin{equation}\label{eq:OBS}\tag{$\mathrm{OBS}$}
	\min\big\{\,\mathrm{L}u,\ u-\bar{u}\,\big\}\,=\,0\,,\qquad
	u(0,\cdot)\,=\,u_0(\cdot)\,.
\end{equation}

Given \eqref{eq:OBS_setting}-\eqref{eq:ellipticity}, the existence of  
viscosity solutions to \OBS$(\sigma,u_0,\bar{u})$ follows from standard 
results (\citealp[see e.g.][]{ElKarouiPardouxPeng:1997}). We then define
\begin{align}\label{eq:defn_D}
	D\ =\ \big\{(t,x)\in\R_+\times\R:\ u(t,x)\,>\,\bar{u}(x)\big\}.
\end{align}
Obviously, $D$ is an open set since $u$ and $\bar{u}$ are continuous.
	
Moreover, if $D=\R_+\times\R$, then $u> \bar{u}$ everywhere. However, 
since $Lu=0$ on $D=\R_+\times\R$, $u(t,x)=-\E^\nu|x-X_t|=\poten{\nu}(x)
-\mathrm{L}^x_t\searrow-\infty$ as $t\to\infty$, which violates the fact 
that $u>\bar{u}$ everywhere. Therefore, we have that $D\subsetneq 
\R_+\times\R$.
	
In this section, we will see that $D$ is a Root's barrier such that
the first hitting time $\tau_D=\inf\big\{\,t>0: (t,X_t)\notin D\,\big\}$ 
is a solution for \SEP$(\sigma,\nu,\mu)$.

The key observation is that the solution $u(t,x)$ has an interpretation in 
terms of an optimal stopping problem 
(\citealp[see][Sect.\,3.4.9]{BensoussanLions:1982}):
\begin{equation*}
	u(t,x) = \sup\nolimits_{\theta\leq t}J_{t,x}(\theta),
		~~\text{where }~J_{t,x}(\theta):=\E^x\big[u_0(Y_\theta)
		+(\bar{u}-u_0)(Y_\theta)\ind_{\theta<t}\big].
\end{equation*}
Here, $Y$ is an independent copy of $X$, but runs backward from $(t,x)$.
Moreover, according to \citet[Rmk\,4.4]{CoxWang:2013a},
\begin{align*}
	u(t,x) = J_{t,x}(\theta_t),\ \ \text{ where }\
		\theta_t=\inf\big\{r\geq0:\ (t-r,Y_r)\notin D\big\}\wedge t.
\end{align*}
Using this result we firstly verify that the open set $D$ 
is a Root's barrier. 
	
\begin{lem}\label{lem:prepare}
	Suppose that \eqref{eq:OBS_setting}-\eqref{eq:ellipticity} hold, then
	$u(t,x)$ is non-increasing in $t$ and $D$ is a Root's barrier.
\end{lem}

\begin{proof}
	For any fixed $(t,x)$, and a stopping time $\theta\leq t$ and a 
	deterministic time $s\leq t$, 
	\begin{equation*}
		\begin{aligned}
			J_{t,x}(\theta)\,&=\,\E^x\big[u_0(Y_{\theta})
				+(\bar{u}-u_0)(Y_\theta)\ind_{\theta<s}
				+(\bar{u}-u_0)(Y_\theta)\ind_{s\leq\theta<t}\big]\\[6pt]
			&\leq\,\E^x\big[u_0(Y_{s\wedge\theta})
				+(\bar{u}-u_0)(Y_{s\wedge\theta})\ind_{s\wedge\theta<s}\big]
				+\E^x\big[u_0(Y_{\theta})-u_0(Y_{s\wedge\theta})\big]\\[6pt]
			&=\,J_{s,x}(s\wedge\theta)+\E^x\big[u_0(Y_{\theta})
				-u_0(Y_{s\wedge\theta})\big]
		\end{aligned}
	\end{equation*}
	where the inequality holds because $\bar{u}\leq u_0$. Then 
	$J_{t,x}(\theta)\leq J_{s,x}(s\wedge\theta)$ by Jensen's inequality
	since $u_0$ is concave. It follows that 
	\begin{equation*}
		u(t,x)\ =\ \sup\nolimits_{\theta\leq t}J_{t,x}(\theta)
			\ \leq\ \sup\nolimits_{\theta\leq s}J_{s,x}(\theta)\ =\ u(s,x).
	\end{equation*}
	Thus, $u(t,x)$ is non-increasing in $t$. It follows that $D$ is a 
	Root's barrier. 
\end{proof}

The non-increase of $u$ in time also can be found in
\citet[Cor.\,1]{GassiatOberhauserReis:2015}, and they proved the result 
using PDE theory. The proof we present here is independently derived via
the connection between optimal stopping problems and obstacle problems.

Next we will interpret the viscosity solution $u(t,x)$ in a probabilistic 
viewpoint.
	
\begin{lem}\label{lem:u_distribution}
	Suppose that \eqref{eq:OBS_setting}-\eqref{eq:ellipticity} hold, 
	then there exists some probability distribution $\mu_t$ such that
	$u(t,\cdot)=\poten{\mu_t}$ for all $t\geq 0$.
\end{lem}

\begin{proof}
	Firstly, the concavity of $u$ in space easily follows from
	the non-increase of $u$ in time and \eqref{eq:ellipticity}.
		
	Noting that $|(u_0)'_-|\leq 1$ and the Radon measure 
	$u_0''(\dif x)=-2\nu(\dif x)$, we have, by It\^o-Tanaka formula,
	\begin{align*}
		0\ \leq\ u_0(x)-u(t,x)\ &\leq\ u_0(x)-J_{t,x}(t)
			\ =\ u_0(x)-\E^x\big[u_0(Y_t)\big]\\[6pt]
		&=\ -\E^x\left[\int_0^t (u_0)'_-(Y_s)\dif Y_s+
			\half\int_\R L^a_t\ u_0''(\dif a)
			\right]\ =\ \int_\R\E^x\big[L^a_t\big]\nu(\dif a).
	\end{align*}
	Denote the transition density of $Y$ by $p^Y_t(x,y)$. By the symmetry	
	property of density (c.f. \citealp[p.149]{ItoMcKean:1974}; 
	\citealp[Thm.\,2.2]{EkstromTysk:2011}),
	\begin{align*}
		\E^x\big[L_t^a\big]\ =\ \int_0^t\sigma^2(a)p^Y_s(x,a)\dif s
			\ =\ \int_0^t\sigma^2(x)p^Y_s(a,x)\dif s
			\ =\ \E^a\big[L^x_t\big].
	\end{align*}
	It then follows from \citet[Lem.\,2.2]{Chacon:1977} that, as 
	$|x|\to\infty$,
	\begin{align*}
		u_0(x)-u(t,x)\,&
			\leq\,\E^\nu\big[L^x_t\big]
			\,=\,\E^\nu\big[|x-X_t|-|x-X_0|\big]\,\longrightarrow\,0.
	\end{align*}
	Thus, we conclude that there exists some probability distribution,
	denoted by $\mu_t$, such that $u(t,\cdot)\,=\,\poten{\mu_t}
	\,\leq\,\poten{\nu}$ (c.f. \citealp[Lem.\,2.3.1]{J.Wang:Thesis}).
\end{proof}
	
Noting that $u(t,x)$ is non-increasing in $t$ and bounded below by 
$\bar{u}(x)$, we can define $\widehat{u}(x)=\lim_{t\to\infty}u(t,x)$ for all 
$x\in\R$. According to \citet[Lem.\,2.5\,\&\,2.6]{Chacon:1977}, there exists 
some constant $C_L$ and a measure $\widehat{\mu}$ defined on $\R$ such that
\begin{align}\label{eq:defn_hat_mu}
	\mu_t\,\Longrightarrow\,\widehat{\mu}\quad\text{and}\quad  
		\widehat{u}(x)\,=\,\poten{\widehat{\mu}}(x)-C_L,
		\quad\forall\ x\in\R\,.
\end{align}
We also define
\begin{align*}
	\widehat{D}\,=\,\big\{(t,x):u(t,x)>\widehat{u}(x)\big\}
		\quad\text{and}\quad\widehat{\tau}\,=\,
		\inf\big\{t>0:(t,X_t)\notin\widehat{D}\big\}.
\end{align*}
Obviously $\widehat{D}\subset D$ and $\widehat{\tau}\leq\tau_D$. 
Moreover, we have the following result.
	
\begin{lem}\label{lem:Q_empty}
	Suppose that \eqref{eq:OBS_setting}-\eqref{eq:ellipticity} hold,
	then $X_{\widehat{\tau}}\sim\widehat{\mu}$.
\end{lem}

\begin{proof}
	For some fixed time $t>0$, one can easily check that $u(\cdot\wedge t,
	\cdot)$ is the viscosity solution of \OBS$(\sigma,\poten{\nu},
	u(t,\cdot))$. Define
	\begin{align*}
		D_t:=\big\{(s,x):u(s\wedge t,x)>u(t,x)\big\},\quad
			\tau_t=\inf\big\{s>0:(s,X_s)\notin D_t\big\}\,.
	\end{align*}
	Then by \thmref{thm:Gassiat}, 
	\begin{align}\label{eq:poten_mu_t}
		u(t,x)\,=\,-\E^\nu\big|x-X_{\tau_t}\big|\quad\text{and}
		\quad X_{\tau_t}\sim \mu_t
	\end{align}
	Since $u$ is non-increasing in time, it is easy to check that 
	$\{D_t\}_{t>0}$ is a non-decreasing sequence of open sets. Further, 
	since $u(t,x)\geq \widehat{u}(x)$, one can check that 
	$D_t\subset \widehat{D}$.
	
	Conversely, for any $(t,x)\in \widehat{D}$, $u(t,x)>\widehat{u}(x)$. 
	Since $\lim_{s\to\infty}u(s,x)=\widehat{u}(x)$, there must be some
	$T>t$ such that $u(t,x)>u(T,x)>\widehat{u}(x)$, i.e. $(t,x)\in D_T$. 
	As conclusion, we have that $D_t\nearrow\widehat{D}$.
	It follows from \eqref{eq:defn_hat_mu} and \eqref{eq:poten_mu_t} that 
	\begin{align}\label{eq:tau_limit}
		\tau_t\,\nearrow\,\widehat{\tau}~~\text{as}~t\to\infty,
			~~\pb^\nu\text{-a.s.},\quad\text{and hence,}
			\quad X_{\widehat{\tau}}=\lim\nolimits_{t\to\infty}
			X_{\tau_t}\sim\widehat{\mu}.
	\end{align}
\end{proof}
	
We then can present the main result of this section, which connects 
Skorokhod embedding problems to obstacle problems when the convex ordering 
condition \eqref{eq:regular} fails.
	
\begin{thm}\label{thm:mainResult}
	Given \eqref{eq:OBS_setting}-\eqref{eq:ellipticity}, let $\tau_D$ be 
	the stopping time defined in \eqref{eq:defn_D}. Then $\tau_D$ is a 
	Root's solution to \SEP$(\sigma,\nu,\mu)$. Moreover, $u(t,x)
	=-\E^{\nu}\big|x-X_{t\wedge\tau_D}\big|$ and $u(t,x)\searrow 
	\bar{u}(x)$ as $t\to\infty$ for all $x\in\R$.
\end{thm}

\begin{proof}
	First of all, we define $F:=\big\{x\in\R:R(x)<+\infty\big\}$, 
	$\widehat{F}:=\big\{x\in\R:\widehat{R}(x)<+\infty\big\}$ where $R$ and 
	$\widehat{R}$ denote the barrier functions of $D$ and $\widehat{D}$ 
	respectively. Because $\widehat{D}\subset D$, we have that $\widehat{R}
	\leq R$ and $\widehat{F}\supset F$. In addition, $F$ is non-empty 
	since $D\subsetneq\R_+\times\R$ as mentioned before, and hence both the
	stopping times $\tau_D$ and $\widehat{\tau}$ are non-trivial, i.e.
	finite almost surely (\citealp[Prop.\,1]{Loynes:1970}).
	
	For any $x\in \widehat{F}\cap F^\complement$, we have that 
	$\widehat{u}(x)=u(\widehat{R}(x),x)>\bar{u}(x)$. By the continuity of 
	$\widehat{u}$ and $\bar{u}$, there exists $\eps$ such that $u_0(y)
	\geq\widehat{u}(y)>\bar{u}(y)$ for all $y\in(x-\eps,x+\eps)$, and then
	$u(t,y)>\bar{u}(y)$ for all $t> 0$. It follows that
	$(0,+\infty)\times(x-\eps,x+\eps)\subset D$. Since $\mathrm{L}u=0$ on 
	$D$, the process $\{u(t-r,Y_r)\}$ is a martingale up to the hitting 
	time $H_{x\pm\eps}$ under $\pb^x$. Therefore, for $t>0$, since $u$ is 
	non-increasing in $t$,
	\begin{equation*}
		u(2t,x)\ =\ \E^x\big[u\big(2t-t\wedge H_{x\pm\eps},
			Y_{t\wedge H_{x\pm\eps}}\big)\big]\ \leq\
			\E^x\big[u(t,Y_{t\wedge H_{x\pm\eps}})\big].
	\end{equation*}
	Let $t\to\infty$, since $u(t,x)\searrow \widehat{u}(x)$, 
	by the concavity of $\widehat{u}$ and Fatou's Lemma,
	\begin{equation*}
		\widehat{u}(x)\ \leq\ \E^x\big[\lim\nolimits_{t\to\infty}
			u(t,Y_{t\wedge H_{x\pm\eps}})\big]\ =\ 
			\E^x\big[\widehat{u}(Y_{H_{x\pm\eps}})\big]
			\ \leq\ \widehat{u}(x).
	\end{equation*}
	Hence $\E^x\big[\widehat{u}(Y_{H_{x\pm\eps}})\big]=\widehat{u}(x)$. 
	Since $u_0>\bar{u}$ on $(x-\eps,x+\eps)$, the process $Y$ is 
	non-degenerate before $H_{x\pm\eps}$ by \eqref{eq:ellipticity}, so
	the concave function $\widehat{u}$ is in fact linear on 
	$(x-\eps,x+\eps)$(c.f. \citealp[Prop.\,3.5.1]{Lange:2010}). This 
	implies that $\widehat{\mu}(x-\eps,x+\eps)=0$, and then it follows that 
	$\widehat{\mu}\big(F\big)=\widehat{\mu}\big(\widehat{F}\big)$. Moreover,
	since $\widehat{R}(X_{\widehat{\tau}})\leq\widehat{\tau}<\infty$ 
	almost surely, we have that $\widehat{\mu}(F)=\widehat{\mu}
	(\widehat{F})=1$ by
	\begin{equation*}
		\widehat{\mu}\big(\widehat{F}\big)\ =\ 
			\pb^\nu\big[\widehat{R}(X_{\widehat{\tau}})<\infty\big]
			\ \geq\ \pb^\nu\big[\widehat{\tau}<\infty\big]\ =\ 1\,.
	\end{equation*}
	Same argument implies $\mu(F)=1$, and then we have that
	\begin{equation}\label{eq:mu_support}
		\mu(F)\ =\ \widehat{\mu}(F)\ =\ 1,\qquad
			\mu(F^\complement)\ =\ \widehat{\mu}(F^\complement)\ =\ 0\,.
	\end{equation}
		
	For any $x\in F$, $u(t,x)=u(R(x),x)=\bar{u}(x)$ for all $t\geq R(x)$,
	and then $\widehat{u}(x)=\lim_{t\to\infty}u(t,x)=\bar{u}(x)$. Hence
	we have that, by the continuity of $\bar{u}$ and $\widehat{u}$,
	\begin{equation*}
		\bar{u}\ =\ \poten{\mu}-C\ \leq\ \poten{\widehat{\mu}}-C_L
			\ =\ \widehat{u}\quad
			\text{ on }\,\R,\quad~~\text{with\ \,``=''\ \,on\,\ $\cl{F}$},
	\end{equation*}
	where $\cl{F}$ denotes the closure of $F$.
	
	Define $x_\ast=\inf F$, $x^\ast=\sup F$. Then we have that 
	\begin{equation*}
		\bar{u}(x_\ast)\,=\,\widehat{u}(x_\ast)\quad
			\text{ if }\ x_\ast>-\infty;~\quad~
		\bar{u}(x^\ast)\,=\,\widehat{u}(x^\ast)\quad
			\text{ if }\ x^\ast<+\infty.
	\end{equation*}
		
	For any $x$ such that $-\infty<x<x_\ast$, since $\mu(F)
	=\widehat{\mu}(F)=1$, it is easy to compute that $\bar{u}(x)\,=\,
	-(m_\mu+C)+x$ and $\widehat{u}(x)\,=\,-(m_{\widehat{\mu}}+C_L)+x$,
	where the mean values of $\mu$ and $\widehat{\mu}$ are denoted by 
	$m_\mu$ and $m_{\widehat{\mu}}$. Let $x\to x_\ast$. It follows from the
	continuity of potential functions and $\bar{u}(x_\ast)=
	\widehat{u}(x_\ast)$ that $-(m_\mu+C)=-(m_{\widehat{\mu}}+C_L)=
	\bar{u}(x_\ast)-x_\ast$, and then $\bar{u}$ and $\widehat{u}$ agree on 
	$(-\infty,x_\ast)$:
	\begin{align*}
		\bar{u}(x)\,=\,\widehat{u}(x)\,=\,\bar{u}(x_\ast)-(x_\ast-x)
		\quad\text{for all }x<x^\ast.
	\end{align*}
	Similarly, we have that  $\bar{u}$ and $\widehat{u}$ also agree on 
	$(x^\ast,+\infty)$:
	\begin{align*}
		\bar{u}(x)\,=\,\widehat{u}(x)\,=\,\bar{u}(x^\ast)-(x-x^\ast)
		\quad\text{for all }x>x^\ast.
	\end{align*}
		
	For the case where $x\in F^{\complement}$ and there exist $z_1,z_2\in 
	F$ such that $z_1<x<z_2$, denote $z_\ast:=\sup\{y\in F:\,y<x\}$, 
	$z^\ast=\inf\{y\in F:\,y>x\}$. Since $(z_\ast,z^\ast)\subset 
	F^\complement$, we have that $\mu((z_\ast,z^\ast))=0$ and $\widehat{\mu}
	((z_\ast,z^\ast))=0$ by \eqref{eq:mu_support}, which implies that both 
	$\bar{u}$ and $\widehat{u}$ are linear  on $(z_\ast,z^\ast)$. In 
	addition, $\bar{u}(z_\ast)=\widehat{u}(z_\ast)$, $\bar{u}(z^\ast)
	=\widehat{u}(z^\ast)$ because $z_\ast$, $z^\ast\in \cl{F}$. Then we can 
	conclude that $\bar{u}=\widehat{u}$ on $(z_\ast,z^\ast)$:
	\begin{equation*}
		\bar{u}(x)\,=\,\widehat{u}(x)\,=\,\dfrac{z^\ast-x}{z^\ast-z_\ast}
			\,\bar{u}(z_\ast)+\dfrac{x-z_\ast}{z^\ast-z_\ast}
			\,\bar{u}(z^\ast),\quad\text{for all }x\in(z_\ast,z^\ast)
	\end{equation*}
		
	As conclusion, we have that $\bar{u}=\widehat{u}$ on $\R$, and hence 
	$\mu$ and $\widehat{\mu}$ agree on $\R$ and so $C=C_L$. It then follows 
	from \eqref{eq:tau_limit} that
	\begin{align*}
		D=\widehat{D},\quad~~X_{\tau_D}=X_{\widehat{\tau}}\sim\mu,
			\quad~~u(t,x)\searrow\bar{u}(x).
	\end{align*}
		
	At last we show that $u(t,x)=-\E^\nu\big|x-X_{t\wedge\tau_D}\big|$. 
	Fix some $t\geq 0$, for all $T>t$, one can easily verify that 
	$u(\cdot\wedge T,\cdot)$ is the viscosity solution to \OBS$(\sigma,
	\poten{\nu},u(T,\cdot))$, and then, by \thmref{thm:Gassiat} and 
	Tanaka's formula, 
	\begin{align*}
		u(t,x)\ =\ -\E^\nu|x-X_{t\wedge\tau_T}|
			\ =\ -|x|-\E^\nu \big[L^x_{t\wedge\tau_T}\big].
	\end{align*}
	Since $\tau_T\nearrow\widehat{\tau}=\tau_D$ (recall 
	\eqref{eq:tau_limit}), let $T\to\infty$, the desired result follows 
	from the monotone convergence theorem.
\end{proof}
	
\begin{rmk}
	We prove the existence and construction of Root's solution to \SEP\ 
	for the case where $X$ is a time-homogeneous diffusion. Thanks to the 
	work of \citet{GassiatOberhauserReis:2015}, our proof also works if
 	the diffusion coefficient $\sigma:\R_+\times\R\to\R$ satisfies 
 	\eqref{eq:Lip_LG} and \eqref{eq:ellipticity} uniformly in 
	time $t$.
\end{rmk}
	
	
\section{Minimality of Root's embeddings}
\label{sec:minimality}
	
In \secref{sec:construction}, we have shown that for any integrable 
distribution $\nu$ and $\mu$, even if \eqref{eq:regular} fails, we still 
can construct a Root's solution to \SEP$(\sigma,\nu,\mu)$ by solving 
\OBS$(\sigma,\poten{\nu},\poten{\mu}-C)$. It also turns out that there exist 
infinitely many Root's embeddings for \SEP$(\sigma,\nu,\mu)$ (dependent on 
different choices of $C$ in the boundary condition).
	
For the cases where $\poten{\mu}\leq \poten{\nu}$, one may think that $C=0$ 
is the best choice because such Root's embeddings are UI stopping times.
For the general cases where $\poten{\mu}\nleq \poten{\nu}$,  we have learned
that there is no UI solution to \SEP$(\sigma,\nu,\mu)$. As mentioned
in \secref{sec:introduction}, now we need the embeddings to be minimal in 
the sense of \citet{Monroe:1972a}. In this section, we study the minimality 
of embeddings, and then, we will see how to choose suitable boundary 
condition in the obstacle problems such that the corresponding Root's 
embeddings are minimal.
	
To this end, we firstly recall the following result 
(\citealp[Thm.\,17]{Cox:2008}), which connects the minimality of stopping 
times to potential functions.
	
\begin{thm}\label{thm:CoxMinimality}
	Let $T$ solve \SEP$(\nu,\mu)$ where $\nu,\mu$ are integrable. Define
	\begin{align}
		\label{eq:defn_A&C}
		\begin{aligned}
			\mathcal{A}\ =\ \big\{\,x\in[-\infty,\,+\infty]:
				\ \lim\nolimits_{y\rightarrow x}\big(\poten{\mu}
				-\poten{\nu}\big)(y)\, =\, C^\ast\,\big\},\\[4pt]
			\text{where }\ \ C^\ast\,:=\,\sup\nolimits_{x\in\R}
				\big\{\poten{\mu}(x)-\poten{\nu}(x)\big\},
		\end{aligned}\\[8pt]
		\label{eq:defn_a_pm} 
		a_+\, =\, \sup\big\{x\in\overline{\R}:\,x\in\mathcal{A}\big\}
			\quad \text{and}\quad a_-\, =\, \inf\big\{x
			\in\overline{\R}:\,x\in\mathcal{A}\big\}.
	\end{align}
	Moreover, denote the first hitting times of the set $\mathcal{A}$ and 
	the horizontal level $\gamma$ by $H_\mathcal{A}$ and $H_\gamma$  
	respectively. Then the following statements are equivalent:
		
	\vspace{-3pt}
	\begin{enumerate}[i), leftmargin = 2.4em, labelsep = .9em ]

		\item\label{itm:Cox_i} $T$ is minimal;
			
		\vspace{-6pt}
		\item\label{itm:Cox_iii} $T\leq H_{\mathcal{A}}$ 
			and for all stopping times $S\leq T$, 
			\begin{equation*}
				\E^\nu\big[W_T\big|\mathcal{F}_S\big]
					\leq W_{S}\,\text{ on }\,\{W_0\geq a_-\};\qquad
				\E^\nu\big[W_T\big|\mathcal{F}_S\big]
					\geq W_{S}\,\text{ on }\,\{W_0\leq a_+\};
			\end{equation*}

		\vspace{-6pt}
		\item\label{itm:Cox_v} $T\leq H_{\mathcal{A}}$ and as 
			$\gamma\rightarrow\infty$, 
			\begin{equation*}
				\gamma\,\pb^\nu\big[T>H_{-\gamma},\ W_0\geq a_-\big]
					\longrightarrow 0;\qquad
				\gamma\,\pb^\nu\big[T>H_{+\gamma},\ W_0\leq a_+\big]
					\longrightarrow 0.
			\end{equation*}
	
	\end{enumerate}
	Further, if there exists $a\in\R$ such that $\pb^\nu[T\leq H_a]=1$, 
	then $T$ is minimal. 
\end{thm}
	
The original proof of \thmref{thm:CoxMinimality} does not rely on any 
properties of Brownian motion beyond the strong Markov property and the 
continuity of paths, so this result can be extended to any continuous 
strong Markov processes.
	
Now, let $\tau$ be a solution to \SEP$(\sigma,\nu,\mu)$ (not necessarily be
of Root's type), we denote the potential of the corresponding stopped 
process by
\begin{equation*}
	u(t,x)\ =\ -\E^\nu\big|\,x-X_{t\wedge\tau}\,\big|.
\end{equation*}
We are interested in what will happen to $u(t,x)$ as $t\rightarrow\infty$.
	
If $\tau$ is a UI stopping time, we immediately have that $\lim_{t\to\infty}
u(t,x)=\poten{\mu}$. For non-UI cases, we firstly review the examples 
mentioned in \secref{sec:introduction}.
	
\begin{example}\label{eg:1}
	For some $a>0$, $H_a=\inf\{t>0:W_t=a\}$ is a non-UI solution for
	\SEP$(\delta_0,\delta_a)$. Let $u(t,x)=-\E^{\delta_0}
	|x-W_{t\wedge H_a}|$. One can compute for $x<a$,
	\begin{align*}
		u(t,x)\ =\ x-2x\cdot\Phi\left(\dfrac{x}{\sqrt{t}}\right)
			+2(x-2a)\cdot\Phi\left(\dfrac{x-2a}{\sqrt{t}}\right)
			-2\sqrt{t}\cdot\left[\,\phi\left(\dfrac{x}{\sqrt{t}}\right)
			-\phi\left(\dfrac{x-2a}{\sqrt{t}}\right)\,\right]\\[8pt]
		\longrightarrow\ x-2x\cdot\Phi(0)+(2x-4a)\cdot\Phi(0)
			\ =\ x-2a\ =\ -|x-a|-a,
	\end{align*}
	where $\Phi$ and $\phi$ denote the CDF and PDF of standard normal 
	distribution respectively. For $x\geq a$, we have that $u(t,x)=-x=
	-|x-a|-a$. Therefore, $\lim_{t\to\infty}u(t,x)=\poten{\delta_a}(x)-a$ 
	for all $x\in\R$.
\end{example}
	
\begin{example}\label{eg:2}
	For some $a>0$, $H_0=\inf\{t>0:W_t=0\}$ is a non-UI solution for
	\SEP$((\delta_a+\delta_{-a})/2,\delta_0)$. Then we have that
	\begin{align*}
		u(t,x)\ =\ -\big(|x|+2a\big)&+\big(|x|+a\big)
			\cdot\Phi\left(\dfrac{|x|+a}{\sqrt{t}}\right)\\[6pt]
		&-\big(|x|-a\big)
			\cdot\Phi\left(\dfrac{|x|-a}{\sqrt{t}}\right)
			+\sqrt{t}\cdot\left[\,\phi\left(\dfrac{|x|+a}{\sqrt{t}}\right)
			-\phi\left(\dfrac{|x|-a}{\sqrt{t}}\right)\,\right].
	\end{align*}
	It is easy to verify that $\lim_{t\to\infty}u(t,x)=-|x|-a
	=\poten{\delta_0}(x)-a$. 
\end{example}
	
\begin{figure}[t]
	
	\caption{The evolution of potentials described 
			in \egref{eg:1} and \egref{eg:2}.}
	\label{fig:examples}

	\subfloat[\quad$\nu=\delta_0,\quad\mu=\delta_1$]{
	\begin{overpic}[trim = 30 30 30 30, clip, width=0.475\textwidth]
		{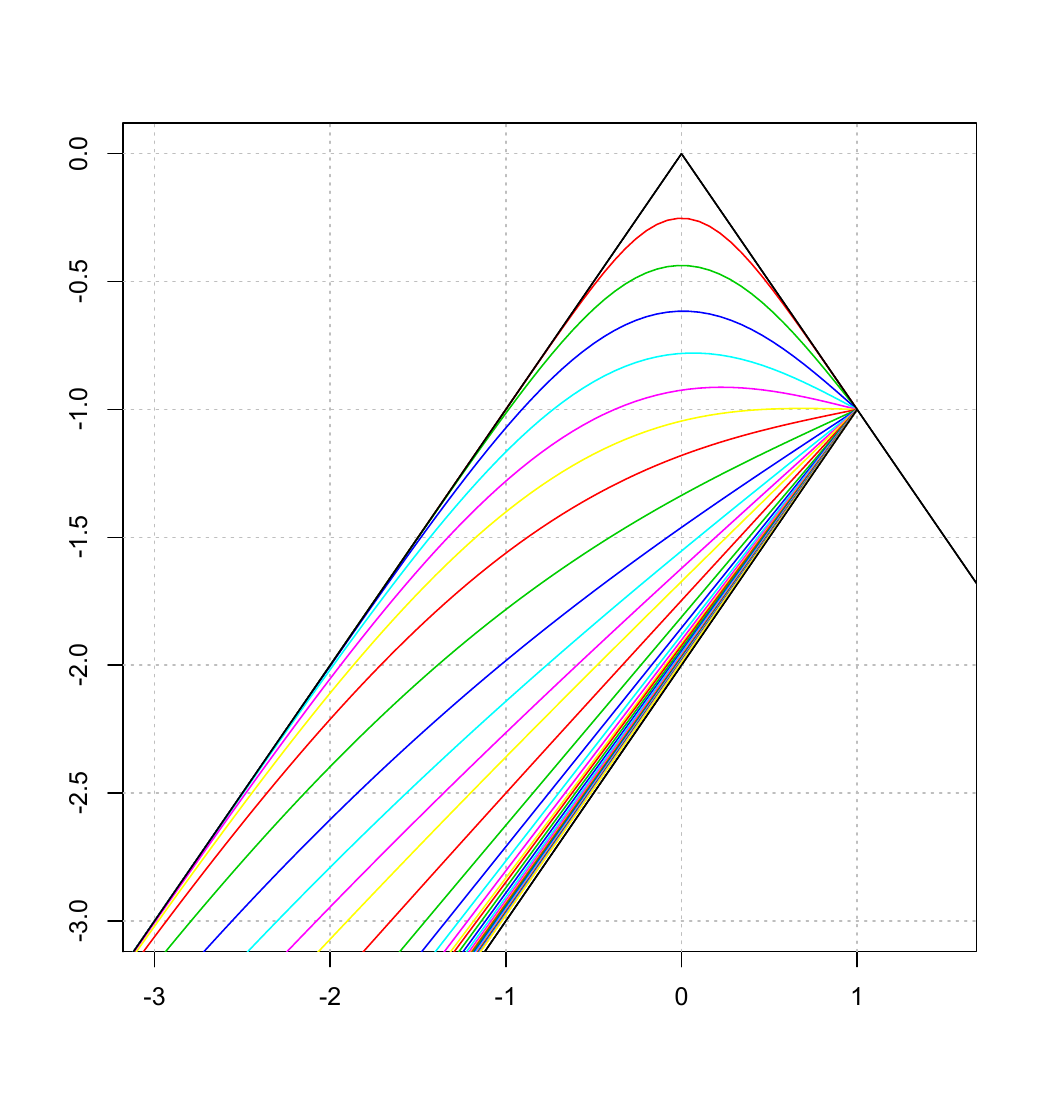}
		\put(24,42){
				\begin{rotate}{56.0}
					$\poten{\nu}=-|x|$
				\end{rotate}}
		\put(51,12){
				\begin{rotate}{56.0}
					$\poten{\mu}-C^\ast=-|x-1|-1$
				\end{rotate}}
	\end{overpic}}
	\hfill
	\subfloat[\quad$\nu=(\delta_1+\delta_{-1})/2,\quad\mu=\delta_0$]{
	\begin{overpic}[trim = 30 30 30 30, clip, width=0.475\textwidth]
		{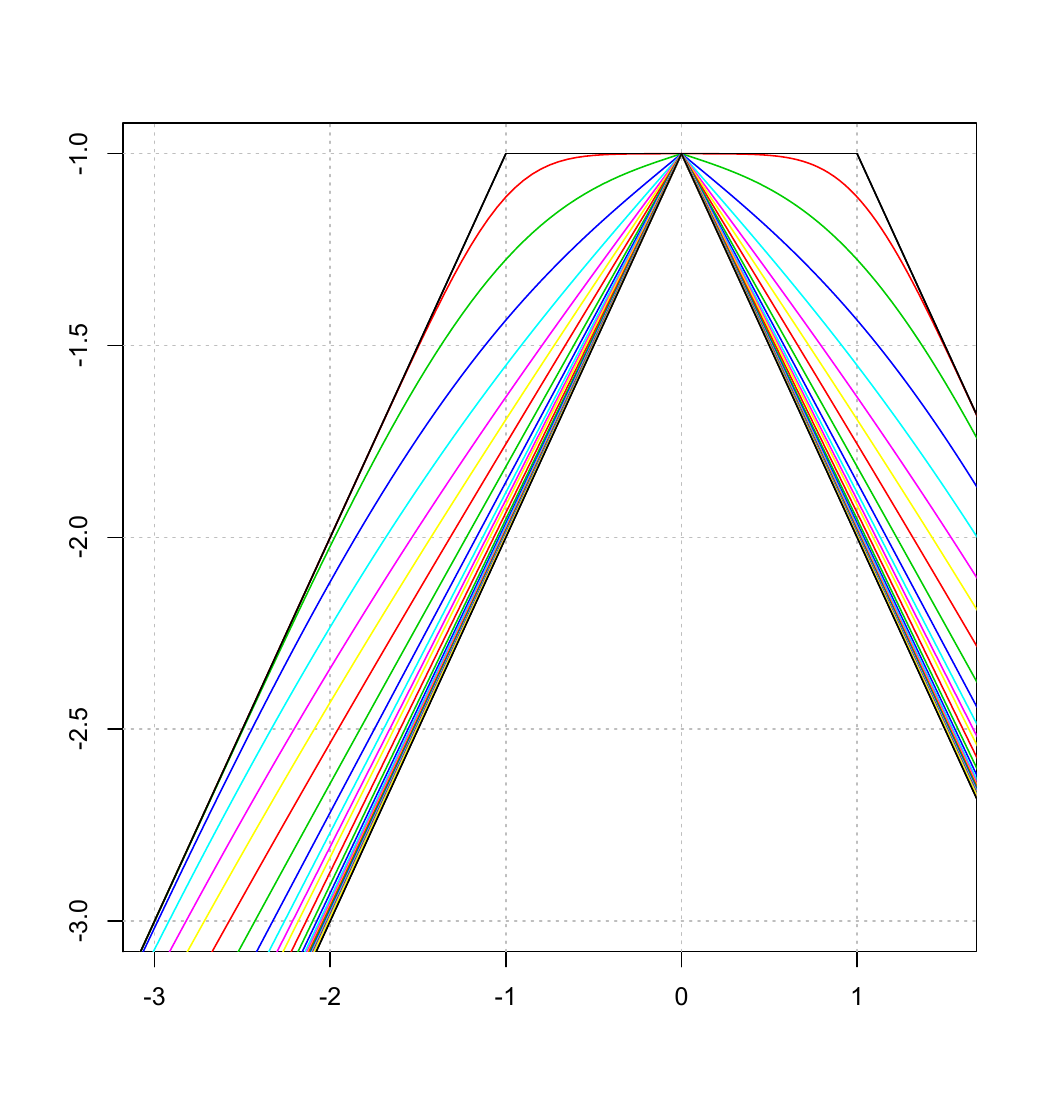}
		\put(15,37){
				\begin{rotate}{66.90}
					$\poten{\nu}=-\max\{|x|,1\}$
				\end{rotate}}
		\put(40,26){
				\begin{rotate}{66.90}
					$\poten{\mu}-C^\ast=-|x|-1$
				\end{rotate}}
	\end{overpic}}
	
	
\end{figure}
	
By the last line of \thmref{thm:CoxMinimality}, both embeddings given in 
the above examples are minimal. Denote the starting and target distributions
by $\nu$ and $\mu$ respectively in these examples, one then can find that
(see \figref{fig:examples})
\begin{align*}
	\lim\nolimits_{t\to\infty}u(t,x)\, =\, \poten{\mu}(x)-C^\ast,\quad 
		\text{where }\, C^\ast = \sup\nolimits_{\R}
		\big\{\poten{\mu}-\poten{\nu}\big\}.
\end{align*}
This result can be extended to general cases as the following lemma.
	
\begin{lem}\label{lem:MinStop_2_MinC}
	Let $\tau$ be a solution to \SEP$(\sigma,\nu,\mu)$ and 
	$C^\ast:=\sup\nolimits_{x\in\R}\big\{\poten{\mu}(x)
	-\poten{\nu}(x)\big\}$. Then
	\begin{equation*}
		\begin{aligned}
				\lim\nolimits_{t\to\infty}u(t,x)\,=\,\poten{\mu}(x)-C_L, 
				\ \ \text{ for all }\,x\in\R\,,\\[6pt]
			\text{where }\ C_L\,=\,C^\ast+\inf\nolimits_{x\in\R}
				\E^\nu\big[L^x_\tau\big].
		\end{aligned}
	\end{equation*}
	In particular, $\,C^\ast=C_L\,$ if $\tau$ is a minimal stopping time.
\end{lem}
\begin{proof}
	Since $t\wedge\tau\rightarrow \tau$, $X_{t\wedge\tau}\rightarrow X_\tau$
	almost surely, and then $\mathcal{L}(X_{t\wedge\tau})\Rightarrow
	\mathcal{L}(X_{\tau})$. By \citet[Lem.\,2.5]{Chacon:1977}, there exists 
	a constant $C_L$ such that
	\begin{align*}
		\lim\nolimits_{t\to\infty}u(t,x)\ =\ \poten{\mu}(x)-C_L,
			\ \ \text{ for all }\,x\in\R.
	\end{align*}
	By martingale property and Tanaka's formula, we have that
	\begin{align*}
		-\infty\ <\ \poten{\mu}(x)-u(t,x)
			\ &=\ \E^\nu\left[\int_{t\wedge\tau}^{\tau}\sgn(x-X_s)\,\dif X_s
			+\left(L^x_{t\wedge\tau}-L^x_{\tau}\right)\right]\\[6pt]
		&=\ \E^\nu\left[\int_{0}^{\tau}\ \sgn(x-X_s)\,\dif X_s
			+\left(L^x_{t\wedge\tau}-L^x_{\tau}\right)\right].
	\end{align*}
	Then, by the monotone convergence theorem,
	\begin{align*}
		C_L\,\equiv\,\poten{\mu}(x)-\lim_{t\to\infty}u(t,x)
			\,=\,\E^\nu\left[\int_0^{\tau}\sgn(x-X_s)\,\dif X_s\right],\ 
			\text{ $\forall\,x\in\R$}.
	\end{align*}
	It follows that, by the definition of $C^\ast$,
	\begin{align*}
		C^\ast\, &=\, \sup\nolimits_{x\in \R}\big\{\poten{\mu}(x)
			-\poten{\nu}(x)\big\}\,=\, \sup\nolimits_{x\in \R}
			\E^\nu\big[|x-X_0|-|x-X_\tau|\big]\\[6pt]
		\,&=\, \sup\nolimits_{x\in \R}\E^\nu\left[\,\int_0^{\tau}
			\!\!\sgn(x-X_s)\dif X_s-L^x_{\tau}\,\right]\,=\,
			C_L-\inf\nolimits_{x\in\R}\E^\nu\big[L^x_{\tau}\big].
	\end{align*}
		
	Now we assume additionally that $\tau$ is a minimal stopping time.
	Consider the following cases dependent on the intersection of $\R$ and 
	$\mathcal{A}$ defined in \eqref{eq:defn_A&C}.

	\vspace*{-2pt}
	\begin{itemize}[leftmargin = 1.56em, labelsep = .8em]
			
		\item \emph{The case where $\mathcal{A}\cap \R\neq \emptyset$.}
			
			We can pick $y\in\mathcal{A}\cap\R$. Since $\tau$ is minimal, 
			by \thmref{thm:CoxMinimality}, we have that
			$\tau\leq H_{\mathcal{A}}\leq H_y<\infty$, $\pb^\nu$-a.s..
			It follows that $\E^\nu\left[L^y_\tau\right]\leq\E^\nu
			\big[L^y_{H_y}\big]=0$. Therefore, $\inf_{x\in \R}\E^\nu
			\left[L^x_{\tau}\right]=0$.
			
		\item \emph{The case where $\mathcal{A}\cap \R= \emptyset$.}
			
			Without loss of generality, we assume that $+\infty\in
			\mathcal{A}$. For any $y\in\R$, denoting $a^+:=\max(a,0)
			=(|a|+a)/2\,$, then
			\begin{equation}\label{eq:+tive}
				\begin{aligned}
					\E^\nu \left[L_{t\wedge\tau}^{y}\right]
						\ &=\ \E^\nu\big|X_{t\wedge\tau}-y\big|
						-\E^\nu\big|X_0-y\big|\\[6pt]
					\ &=\ 2\E^\nu\big[(X_{t\wedge\tau}-y)^+\big]
						-2\E^\nu\big[(X_{0}-y)^+\big].
				\end{aligned}
			\end{equation}
			Since $\tau$ is a minimal stopping time, by Jensen's inequality 
			and \thmref{thm:CoxMinimality}\,\itmref{itm:Cox_iii},
			\begin{equation*}
				(X_{t\wedge\tau}-y)^+\ \leq\ \Big(\E^\nu\Big[X_{\tau}
					-y\,\Big|\mathcal{F}_{t\wedge\tau}\Big]\Big)^+
					\ \leq\ \E^\nu\Big[(X_{\tau}-y)^+\Big|
					\mathcal{F}_{t\wedge\tau}\Big].
			\end{equation*}
			Then the process $\big\{(X_{t\wedge\tau}-y)^+\big\}$ is 
			uniformly integrable because $\mu$ is integrable. Now letting 
			$t$ go to infinity in \eqref{eq:+tive}, we have that
			\begin{align*}
				\E^\nu \left[L_{\tau}^{y}\right]\ =&\ 
					2\E^\nu\big[(X_{\tau}-y)^+\big]-2\E^\nu
					\big[(X_{0}-y)^+\big]\\[6pt]
				=&\ \Big\{\E^\nu|X_{\tau}-y|+\E^\nu[X_{\tau}-y]\Big\}
					-\Big\{\E^\nu|X_{0}-y|+\E^\nu[X_{0}-y]\Big\}\\[6pt]
				=&\ \big[\,\poten{\nu}(y)-\poten{\mu}(y)\,\big]
					+(m_\mu-m_\nu),
			\end{align*}
			where $m_\nu$ and $m_\mu$ denote the mean values of $\nu$ and 
			$\mu$ respectively. On the other hand, by 
			\citet[Lem.\,2.2]{Chacon:1977}, we have that
			\begin{align*}
				\lim_{y\to+\infty}\big[\poten{\nu}(y)+(y-m_\nu)\big]\ =\
				\lim_{y\to+\infty}\big[\poten{\mu}(y)+(y-m_\mu)\big]\ =\ 0,
			\end{align*}
			which implies that $\poten{\mu}(y)-\poten{\nu}(y)
			\longrightarrow m_\mu-m_\nu$ as $y\to+\infty$. 
			We then can conclude that 
			$\E^\nu \left[L_{\tau}^{y}\right]\rightarrow0$ as $y\to+\infty$, 
			and then $\inf_{y\in\R}\E^\nu \left[L_{\tau}^{y}\right]=0$.
			The case where $-\infty\in\mathcal{A}$ is similar.
	
	\end{itemize}
\end{proof}	

We have seen that $C^\ast=C_L$ (or equivalently, $\inf_{x\in\R}\E^\nu
\big[L^x_\tau\big]=0$) is a necessary condition for the minimality.
However, our aim in this section is to show that the Root's embedding 
given by \OBS$(\sigma,\poten{\nu},\poten{\mu}-C^\ast)$ is minimal. To this 
end, next we will see $C^\ast=C_L$ is also a sufficient condition.
	
\begin{thm}\label{thm:minimal}
	Under the same assumptions imposed in \lemref{lem:MinStop_2_MinC}, 
	$\tau$ is a minimal stopping time if and only if $\,\inf_{x\in\R}
	\E^\nu\big[L^x_\tau\big]=0$, or equivalently,
	\begin{align*}
		\lim\nolimits_{t\to\infty}u(t,x)=\poten{\mu}(x)-C^\ast,\quad
		\text{for all}\,~x\in\R.
	\end{align*}
\end{thm}

\begin{proof}
	It has been shown in \lemref{lem:MinStop_2_MinC} that 
	$\lim_{t\rightarrow\infty}u(t,x)=\poten{\mu}(x)-C^\ast$ if $\tau$ is 
	minimal.It only remains to show the ``if'' part. Now we suppose that 
	$\lim_{t\rightarrow\infty}u(t,x)=\poten{\mu}(x)-C^\ast$. Consider the 
	following cases dependent on the intersection of $\R$ and $\mathcal{A}$
	defined in \eqref{eq:defn_A&C}.
		
	\vspace*{-6pt}
	\begin{itemize}[leftmargin = 1.56em, labelsep = .8em ]
			
		\item\emph{The case where $\mathcal{A}\cap \R\neq \emptyset$.}
			
			We can pick  some $y\in\mathcal{A}\cap\R$. Since potential 
			functions are continuous and $u(t,x)\rightarrow\poten{\mu}(x)
			-C^\ast$, we have that $\lim\nolimits_{t\to\infty}u(t,y)
			=\poten{\nu}(y)$. Then by Tanaka's formula and monotone 
			convergence theorem,
			\begin{align*}                           
				\E^\nu\big[ L_{\tau}^{y}\big]\ &=\ \lim_{t\rightarrow\infty}
					\E^\nu\big[ L_{t\wedge\tau}^{y}\big]\ =\ \poten{\nu}(y)
					-\lim_{t\rightarrow\infty}u(t,y)\ =\ 0,     
			\end{align*}                             
			and hence, $L^y_\tau=0$, $\pb^\nu$-a.s.. It follows that
			$\tau\leq H_{y}$ almost surely, and then $\tau$ is a minimal 
			stopping time by the last line of \thmref{thm:CoxMinimality}.
			
		\item \emph{The case where $\mathcal{A}\cap \R= \emptyset$.}
			
			Suppose that $+\infty\in\mathcal{A}$. Since $C^\ast=
			\lim_{y\to+\infty}[\poten{\mu}(y)-\poten{\nu}(y)]= m_\mu-m_\nu$
			(see the proof of \lemref{lem:MinStop_2_MinC}) and $u(t,x)
			\rightarrow\poten{\mu}(x)-C^\ast$, we have that   
			\begin{align*}                          
				2\E^\nu\big[X_{\tau}^+-X_{t\wedge\tau}^+\big]
					\ =\ &\E^\nu\big[|X_{\tau}\big|-|X_{t\wedge\tau}|\big]
					+\E^\nu\big[X_{\tau}-X_{t\wedge\tau}\big]\\[6pt]
				\ =\ &\big[\,u(t,0)-\poten{\mu}(0)\,\big]
					+(m_\mu-m_\nu)\\[6pt]
				\ \longrightarrow\ &(m_\mu-m_\nu)-C^\ast\ =\ 0,
				\qquad \text{as }\ t\rightarrow\infty.
			\end{align*}                                                                             
			Then, by Scheff\'e's Lemma, $\big\{X_{t\wedge\tau}^+\big\}$
			is uniformly integrable. Therefore, as $\gamma\to +\infty$,
			\begin{align*}                           
				\gamma\pb^\nu\big[\tau>H_{\gamma}\big]   
					\ &=\ \gamma\cdot\pb^\nu\big[\tau>H_{\gamma},\,X_0
					\geq\gamma\big]
					+\gamma\cdot\pb^\nu\big[\tau>H_{\gamma},\,
					X_0<\gamma\big]\\[6pt]
				\ &\leq\ \gamma\cdot\nu\big([\gamma,\infty)\big)
					+\gamma\cdot\pb^\nu\big[\tau>H_{\gamma},\,X_0
					<\gamma\big]\\[6pt]
				\ &\leq\ \gamma\cdot\nu\big([\gamma,\infty)\big)+
					\gamma\cdot\pb^\nu\big[X_{H_\gamma \wedge \tau}
					\geq\gamma,\,X_0<\gamma\big]\\[6pt]
				\ &\leq\ \E^\nu\big[X_0;\,X_0\geq\gamma\big]
					+\E^\nu\left[X_{H_\gamma \wedge \tau};\,
					X_{H_\gamma \wedge \tau}\geq\gamma\right]
					\ \longrightarrow\ 0.
			\end{align*}
			Then it follows from 
			\thmref{thm:CoxMinimality}\,\itmref{itm:Cox_v} 
			that $\tau$ is minimal. The case in which 
			$-\infty\in\mathcal{A}$ is similar.
	\end{itemize}
\end{proof}
	                       
Thanks to \thmref{thm:minimal}, we can directly tell if the Root's embedding 
given by \thmref{thm:mainResult} is minimal or not.
	
\begin{thm}\label{thm:minimalResult}
	For integrable probability distributions $\nu$ and $\mu$ on $\R$, 
	assume that $C^\ast=\sup_{x\in\R}\big\{\poten{\mu}(x)-
	\poten{\nu}(x)\big\}$ and $\sigma$ satisfies \eqref{eq:Lip_LG} and 
	\eqref{eq:ellipticity}. Let $u(t,x)$ be the viscosity solution to 
	\OBS$(\sigma,\poten{\nu},\poten{\mu}-C^\ast)$, and $D$ be the set 
	defined in \eqref{eq:defn_D}. Then $\tau_D$ is a minimal solution to 
	\SEP$(\sigma,\nu,\mu)$. Moreover, we have the presentation
	that $u(t,x)=-\E^\nu\big|x-X_{t\wedge\tau_D}\big|$.
\end{thm}

	
\section{Optimality of minimal Root's embeddings}
\label{sec:optimality}
	
As well-known, the UI embedding of Root's type is remarkable because it is 
of minimal residual expectation (\ref{eq:mre}). A natural question now 
arises: can we generalize this optimality result to non-UI Root's 
embeddings? When the stopped process $X^{\tau}$ is not uniformly integrable, 
we cannot expect that $\E^\nu[(\tau-t)^+]$ is finite. Thus, we study the
quantity $\E^\nu[\tau\wedge t]=\E^\nu[\tau-(\tau-t)^+]$ instead. We 
conjecture that the minimal Root's embedding $\tau$ is of \emph{maximal 
principal expectation}, that is, 
\begin{align*}
	\textit{Amongst all minimal solutions of
		 \,\ref{eq:SEP}$\,(\sigma,\nu,\mu)$, the Root's solution}&\\[4pt]
	\textit{maximises 
		$\E^\nu\big[\tau\wedge t\big]$ simultaneously for all $t>0$}&.
\end{align*}
	
For \ref{eq:SEP}$(1,\nu,\mu)$, given that \eqref{eq:regular} holds, this 
statement holds obviously since the minimal Root's solution is UI and then 
is of \ref{eq:mre}.
	
For general cases in which \eqref{eq:regular} fails, we suppose that $\tau$ 
is a minimal embedding for \SEP$(\sigma,\nu,\mu)$, and the stopped potential 
is denoted by $u^\tau(t,x)\,=\,-\E^\nu\big|x-X_{t\wedge\tau}\big|$.
	
It is obvious that $u^\tau(t,x)$ is non-increasing in $t$. According to 
\thmref{thm:minimal}, $u^\tau(t,x)\rightarrow\poten{\mu}(x)-C^\ast$, and 
hence $u^\tau(t,x)\geq \poten{\mu}(x)-C^\ast$ for all $(t,x)\in\R_+\times
\R$. On the other hand, since $u^\tau(t,x)=\poten{\nu}(x)-\E^\nu
\left[L^x_{t\wedge\tau}\right]$, we expect that, in the sense of 
distribution\footnote{
	We have to mention that the argument here is just an intuitive 
	illustration without technique details. We shall refer readers to 
	\citet[Thm.\,1]{GassiatOberhauserReis:2015} for a rigorous proof of 
	$\mathrm{L}u^\tau\geq0$.\label{ft:Gassiat}},
\begin{equation}\label{eq:viscosity_super}
	\begin{aligned}
		u^\tau(t+\delta,x)-u^\tau(t,x)\, &=\, -\!\int_t^{t+\delta}
			\dfrac{\sigma^2(x) \pb
			\left[X_s\in\dif x,s<\tau\right]}{\dif x}\dif s\\[6pt]
		&\geq\, -\!\int_t^{t+\delta}\dfrac{\sigma^2(x)
			\pb\left[X_{s\wedge\tau}\in\dif x\right]}{\dif x}\dif s
			\, =\, \dfrac{\sigma^2(x)}{2}\!\int_t^{t+\delta}
			\!\pdd{}{x}u^\tau(s,x)\dif s.
	\end{aligned}
\end{equation}
It then follows that $u^\tau$ is a viscosity supersolution of \OBS$(\sigma,
\poten{\nu},\poten{\mu}-C^\ast)$, while $u^{\tau_D}$ is a viscosity solution 
of \OBS$(\sigma,\poten{\nu},\poten{\mu}-C^\ast)$. According to 
\citet[Thm.\,5]{GassiatOberhauserReis:2015}, we then have the following 
result as an extension of \citet[Thm.\,3]{GassiatOberhauserReis:2015}.
	
\begin{prop}\label{prop:optimality}
	Assume that $\sigma$ satisfies 
	\eqref{eq:OBS_setting}-\eqref{eq:ellipticity}. Let $\tau_D$ and $\tau$ 
	be minimal solutions to \SEP$(\sigma,\nu,\mu)$, among which $\tau_D$ 
	is of Root's type. Then for any $t\geq 0$, $u^{\tau_D}(t,\cdot)\leq 
	u^\tau(t,\cdot)$ on $\R$, or equivalently, $\mathcal{L}(X_{t\wedge\tau})
	\preceq\mathcal{L}(X_{t\wedge\tau_D})$ in convex order.
\end{prop}
	
Beyond the optimality result above, we are also interested in deriving a 
pathwise inequality which encodes the maximal principal expectation in the 
sense that, for a non-decreasing \emph{concave} function $F:\R_+\to\R_+$ 
with $F(0)=0$, we can find a supermartingale $G_t$ and a function $H(x)$ 
such that $F(t) \leq G_t+H(X_t)$, and such that the equality holds when 
$t=\tau_D$, and $\{G_{t\wedge\tau_D}\}$ is a martingale.
	
Different from our work in \citet{CoxWang:2013a}, here we do not assume 
distributions are in convex order (or equivalently, the embeddings are UI
stopping times), so we cannot take limit on $G_{t\wedge\tau}$ as before. 
Instead, in the following proof, we will see the limit of 
$G_{t\wedge\tau_D}-G_{t\wedge\tau}$ does exist even if the embeddings are 
not UI. 
	 
\begin{thm}\label{thm:path_optimality}
	Suppose that $\sigma(\cdot)$ satisfies 
	\eqref{eq:OBS_setting}-\eqref{eq:ellipticity} and
	\begin{equation}\label{eq:sigma_condition}
		\E^\nu\left[\int_0^{X_0}\!\dif y\!\int_0^y
			\dfrac{\dif z}{\sigma^2(z)}+\int_0^T\!\!\left(
			\int_0^{X_s}\dfrac{\sigma(X_s)}{\sigma(z)}\dif z\right)^{\!\!2}
			\dif s\right]<\infty,\quad\text{for all}~T>0.
	\end{equation}
	Let $\tau$ and $\tau_D$ be two minimal solutions to \SEP$(\sigma,\nu,
	\mu)$, among which $\tau_D$ is a Root's stopping time. Then $\E^\nu
	\big[F(\tau_D)\big]\geq\E^\nu\big[F(\tau)\big]$ for all non-decreasing 
	concave function $F$.
\end{thm}
	
\begin{proof}	
	Without loss of generality, we always assume that $F(0)=0$. Let
	$f(t)=F'_+(t)$ be the right derivative of $F$. Furthermore we may 
	assume that $f$ is bounded and vanishes after some time:
	\begin{align}\label{eq:vanish_f}
		\exists\,N>0,\quad\text{s.t.}\quad\text{(i).}~~f(0)\leq N;
			\quad\text{(ii).}~~\forall\,t>N,~~f(t)\,=\,0.
	\end{align}
	Define $M(t,x)\,:=\,\E^{(t,x)}\big[f(\tau_D)\big]$. Since $f$ is 
	non-increasing, $M(t,x)\leq f(t)$ with equality for $t\geq R(x)$, where 
	$R$ denotes the boundary function of $D$. Hence
	\begin{align}\label{eq:pathwise}
		\begin{aligned}
			\int_0^t\!M(s,x)\dif s+&\int_{0}^{R(x)}
				\!\big[ f(s)-M(s,x) \big]\dif s\\[6pt]
			=\ &\int_0^{R(x)}\!f(s)\dif s-\int^{R(x)}_t
				\!M(s,x)\dif s\ \geq\ F(t),
				\quad\text{with ``='' if}~t\geq R(x).
		\end{aligned}
	\end{align}
	Similar as in our previous work (\citealp[Sect.\,5]{CoxWang:2013a}), 
	we define
	\begin{align*}
		G(t,x):=\int_0^t\!M(s,x)\dif s-Z(x),\quad\text{where}~~
		Z(x)\,:=\,\int_0^x\!\dif y\!\int_0^y\dfrac{2M(0,z)}{\sigma^2(z)}
		\dif z.
	\end{align*}
	Since $M(0,z)\leq f(0)\leq N$, it follows from 
	\eqref{eq:sigma_condition} that
	\begin{equation*}
		\E^\nu\left[Z(X_0)+\int_0^T Z'(X_s)^2\sigma(X_s)^2\dif s\right]
			<\infty,\quad\text{for all }T>0.
	\end{equation*}
	The integrability of $Z(X_t)$ then follows from
	\begin{align*}
		\E^\nu\big[Z(X_{t}) \big]\,=\,\E^\nu\left[Z(X_0)+\int_0^{t} 
			Z'(X_s)\sigma(X_s)\dif W_s
			+\int_0^{t} M(0,X_s)\dif s\right]&\\[6pt]
		\,\leq\ \E^\nu\big[Z(X_0)\big]+N t\,&\,<\ \infty,
	\end{align*}
	and the integrability of $G(t,X_t)$ follows from 
	$\int_0^tM(s,x)\dif s<N^2$ in addition.

	Further, note that we do not need the assumption that $\tau_D$ is UI 
	in the proof of \citet[Lem.\,5.2]{CoxWang:2013a}. Thus, similarly 
	(with a simple modification because the function $f$ is non-increasing 
	here), one can find that $\big\{G(t,X_{t})\big\}$ is a
	$\pb^{\nu}$-supermartingale and a $\pb^{\nu}$-martingale on $\big[0,
	\tau_D\big]$. It then follows from the definition of $G$ that
	\begin{equation}\label{eq:M_tau_ineq}
		\begin{aligned}
			\E^\nu\left[\int_0^{t\wedge\tau_D}
				\!M(s,X_{t\wedge\tau_D})\dif s-\int_0^{t\wedge\tau}
				\!M(s,X_{t\wedge\tau})\dif s\right]\ 
				\geq\ E^\nu\big[Z(X_{t\wedge\tau_D})
				-Z(X_{t\wedge\tau})\big].
		\end{aligned}
	\end{equation}
	Since $Z(0)=Z'(0)=0$, it follows from integration by parts that
	\begin{equation*}
		\begin{aligned}
			Z(x)\,&=\,\int_{\R_+}(x-y)^+\,Z''(y)\dif y
				+\int_{\R_-}(y-x)^+\,Z''(y)\dif y\\[6pt]
			\,&=\,\int_{\R_+}\dfrac{|x-y|+(x-y)}{2}\,Z''(y)\dif y
				+\int_{\R_-}\dfrac{|y-x|+(y-x)}{2}\,Z''(y)\dif y.
		\end{aligned}
	\end{equation*}
	Then for any probability distribution $\lambda$ with mean $m_\lambda$, 
	we have that (assume integrability)
	\begin{equation*}
		\E^{Y\sim\lambda}[Z(Y)]\ =\ \int_{\R_+}
			\!\dfrac{m_\lambda-y-\poten{\lambda}(y)}{2}\,Z''(y)\dif y
			+\!\int_{\R_-}\!\dfrac{y-m_\lambda-\poten{\lambda}(y)}{2}
			\,Z''(y)\dif y.
	\end{equation*}
	Thus, by the fact that $\E^{\nu}\big[X_{t\wedge\tau_D}\big]
	=\E^{\nu}\big[X_{t\wedge\tau_D}\big]=\E^{\nu}\big[X_{0}\big]$, we 
	deduce that
	\begin{equation}\label{eq:Z_expectation}
		\E^\nu\big[Z(X_{t\wedge\tau_D})-Z(X_{t\wedge\tau})\big]\ =\ 
			\int_\R\dfrac{u^\tau(t,y)-u^{\tau_D}(t,y)}{2}Z''(y)\dif y.
	\end{equation}
	By minimality, $u^{\theta}(t,\cdot)\searrow \poten{\mu}-C^\ast$
	for $\theta =\tau_D,\tau$, then
	\begin{equation*}
		\big(\poten{\mu}-C^\ast\big)-u^{\tau_D}\ \leq\ 
			u^\tau-u^{\tau_D}\ \leq\ u^\tau-\big(\poten{\mu}-C^\ast\big).
	\end{equation*}
	Thus it follows from monotone convergence and squeeze theorem that
	\begin{equation*}
		\lim\nolimits_{t\to\infty}\E^\nu\big[Z(X_{t\wedge\tau_D})
			-Z(X_{t\wedge\tau})\big]\ =\ 0.
	\end{equation*}	
	Together with \eqref{eq:M_tau_ineq} and the fact that $\int_0^t M(s,x)
	\dif s\leq  N^2$ for all $(t,x)$ (by \eqref{eq:vanish_f}), it follows 
	from dominated convergence theorem that
	\begin{equation}\label{eq:M_compare}
		\E^\nu\left[\int_0^{\tau_D}M(s,X_{\tau_D})\dif s
			-\int_0^{\tau}M(s,X_{\tau})\dif s\right]\,\geq\,0.
	\end{equation}
		
	On the other hand, since $\tau_D\geq R(X_{\tau_D})$, the inequality 
	\eqref{eq:pathwise} implies that
	\begin{align*}
		F(\tau_D)-F(\tau)\,\geq\,&\int_0^{\tau_D}M(s,X_{\tau_D})\dif s
			-\int_0^{\tau}M(s,X_{\tau})\dif s\\[6pt]
		+&\int_0^{R(X_{\tau_D})}\big[f(s)-M(s,X_{\tau_D})\big]\dif s
			-\int_0^{R(X_{\tau})}\big[f(s)-M(s,X_{\tau})\big]\dif s.
	\end{align*}
	Because $\mathcal{L}(X_{\tau_D})=\mathcal{L}(X_{\tau})$ and all the 
	terms above are integrable by \eqref{eq:vanish_f}, 
	\begin{equation*}
		\E^\nu\big[F(\tau_D)-F(\tau)\big]\ \geq\ 
			\E^\nu\left[\int_0^{\tau_D} M(s,X_{\tau_D})\dif s
			-\int_0^{\tau} M(s,X_{\tau})\dif s\right]\ \geq\ 0.
	\end{equation*}
		
	To observe that the result still holds when \eqref{eq:vanish_f} does 
	not hold, we define $F_N(t)=\min\{Nt,F(t\wedge N)\}$. Then $F_N$ is 
	non-decreasing, concave function satisfying  \eqref{eq:vanish_f}. Hence  
	$\E^\nu[F_N(\tau_D)]\geq\E^\nu[F_N(\tau)]$. Then it follows from the 
	monotone convergence theorem that
	\begin{equation*}
		\E^\nu[F(\tau_D)]\ =\ \lim_{N\to\infty}\E^\nu[F_N(\tau_D)]\
			\geq\ \lim_{N\to\infty}\E^\nu[F_N(\tau)]\ =\ \E^\nu[F(\tau)].
	\end{equation*}	
\end{proof}
	
\begin{rmk}\label{rmk:unsufficient_comparison_1}
	One may find that $\E^\nu\big[Z(X_{t\wedge\tau_D})
	-Z(X_{t\wedge\tau})\big]\geq 0$ because $Z$ is a convex function and 
	$\mathcal{L}(X_{t\wedge\tau})\preceq\mathcal{L}(X_{t\wedge\tau_D})$ in 
	convex order (\propref{prop:optimality}), and then \eqref{eq:M_compare}
	holds by dominated convergence. However, here we show 
	\eqref{eq:M_compare} by showing that $\E^\nu\big[Z(X_{t\wedge\tau_D})
	-Z(X_{t\wedge\tau})\big]$ vanishes as $t$ goes to infinity. The chief 
	reason we adopt such a proof is that the comparison between viscosity 
	(super-)solutions is not sufficient to show the optimality of Root's
	solutions to multi-marginal \mSEP, while the proof presented here still 
	works under such cases (see \secref{sec:multi_optimality}).
\end{rmk}
	

\section{Multi-marginal Skorokhod embedding problem}
\label{sec:multi_margin}
	
In this section, we will extend our results to multi-marginal Skorokhod 
embedding problems. Thanks to a very recent paper, 
\citet{CoxOblojTouzi:2018}, and the arguments presented in previous sections 
of this work, it is not difficult to construct Root's embeddings to such 
multi-marginal embedding problems.
	
\subsection{Construction of Root's embeddings to multi-marginal SEP}
\label{sec:multi_construction}

\citet{CoxOblojTouzi:2018} study the long-standing question of a 
multi-marginal Skorokhod embedding problem \mSEP$(\sigma,\mu_0,\mmu)$ where 
$\mmu$ is a sequence of integrable probability measures $\mu_1,\mu_2,\cdots,
\mu_n$:
\begin{equation}\tag{$\mathrm{SEP}$}\label{eq:mSEP}
	\begin{aligned}
		\textit{Given}~~X_0\sim\mu_0,~~&\textit{to find stopping times}
			~~\tau_1\,\leq\,\tau_2\,\leq\,\cdots\,\leq\,\tau_n,\\[4pt]
		&\textit{such that}~~X_{\tau_1}\sim\mu_1,
			\ X_{\tau_2}\sim\mu_2, \cdots,\ X_{\tau_n}\sim\mu_n.
	\end{aligned}
\end{equation}
	
Given that $\mu_0$ and $\mmu=\{\mu_k\}_{k=1,\cdots,n}$ is of convex 
ordering:
\begin{equation}\label{eq:mregular}
	\poten{\mu_0}(x)\,\geq\,\poten{\mu_1}(x)\,\geq\,\cdots\,\geq\,
		\poten{\mu_{n-1}}(x)\,\geq\,\poten{\mu_{n}}(x)\,>\,-\infty,
		\quad\text{for all}~~x\in\R,
\end{equation}
we consider the following iterated optimal stopping problems: 
\begin{align*}
	u_0(t,x)=\poten{\mu_0}(x),\quad u_k(t,x)=\sup\nolimits_{\theta\leq t}
		J^k_{t,x}(\theta),\quad\text{for }~k=1,2,\cdots,n\,,\\[6pt]
	\text{where}~~J^k_{t,x}(\theta):=\E^x\big[u_{k-1}(t-\theta,Y_\theta)
		+(\poten{\mu_{k}}\!-\poten{\mu_{k-1}})(Y_\theta)
		\ind_{\theta<t}\big].
\end{align*}
Using the solutions $\{u_k\}$, one can define 
\begin{align*}
	\text{for}~~k=1,2,\cdots,n,\quad\tau_0:=0,\quad\tau_k:=\inf
		\big\{t>\tau_{k-1}:(t,X_t)\notin D_{k}\big\},\\[6pt]
	\text{where}~~D_{k}=\big\{(t,x):u_{k}(t,x)-u_{k-1}(t,x)>
		\poten{\mu_k}(x)-\poten{\mu_{k-1}}(x)\big\},
\end{align*}
and then, \citet[Thm.\,3.1]{CoxOblojTouzi:2018} say that, 
\begin{align*}
	\big\{D_k\big\}~\text{are Root's barriers},\quad\big\{\tau_k\big\}
		~\text{is a UI solution to \SEP$(\sigma,\mu_0,\mmu)$},\\[6pt]
	\text{moreover, }\ u_k(t,x)\,=\,-\E^{\mu_0}\big|x-X_{t\wedge\tau_k}\big|
		~~\text{for all}~~k=1,2,\cdots,n.
\end{align*}

Inspired by the result, we are going to consider the multi-marginal \mSEP\ 
when the convex ordering \eqref{eq:mregular} fails.
	
\begin{example}\label{eg:3}
	We shall begin with a simple example \SEP$(\delta_0,
	\{\delta_1,\delta_{-1}\})$. Obviously, a solution of this problem is 
	given by $\tau_1=H_1$, $\tau_2=\inf\{t\geq H_1:X_t= -1\}$. Same as in
	\secref{sec:construction}, we are interested in 
	the limit of $u_j(t,x):=-\E^0|x-W_{t\wedge\tau_j}|$ as $t\to\infty$. 
	For $j=1$, according to \thmref{thm:minimalResult} (or \egref{eg:1}), 
	$u_1(t,x)\to\poten{\delta_1}(x)-1$. For $j=2$,
	\begin{align*}
		u_2(t,x)\ =&\ -\E^0\big[|x-W_{t\wedge\tau_2}|\ind_{t\geq\tau_1}\big]
			-\E^0\big[|x-W_{t}|\ind_{t<\tau_1}\big]\\[6pt]
		=&\ u_1(t,x)-\E^0\big[|x-W_{t\wedge\tau_2}|\ind_{t\geq \tau_1}\big]
			+\E^0\big[|x-W_{\tau_1}|\ind_{t\geq \tau_1}\big]\\[6pt]
		=&\ u_1(t,x)-\int_0^t\wti{\E}^0\Big[\Big|x-1
			-\wti{W}_{\!(t-s)\wedge\wti{H}_{-2}}\Big|\Big]
			\,\pb^0[H_1\!\in\! \dif s]+|x-1|\cdot\pb^0[H_1\leq t]\\[6pt]
		\longrightarrow&\ \big(\!-\big|x-1\big|-1\big)+\big(\!-\big|(x-1)
			-(-2)\big|-2\big)+|x-1|\ =\ \poten{\delta_{-1}}(x)-3.
	\end{align*}
	Here we can interpret the constant $C_2=3$ as
	\begin{equation*}
		C_2\,=\,\sup\nolimits_\R\big\{\poten{\delta_{-1}}-(\poten{\delta_1}
			-C_1)\big\},\quad\text{where}~~
			C_1\,=\,\sup\nolimits_\R\big\{\poten{\delta_1}
			-\poten{\delta_0}\big\}\,=\,1.
	\end{equation*}
	
	We also have to mention that $\tau_2$ is not a minimal stopping time.
	In fact, $W_{H_{-1}}=W_{\tau_2}=-1$ and $H_{-1}\leq\tau_2$ almost 
	surely, and $\pb^0[H_{-1}< \tau_2]>0$. We will see later that the 
	sequence $\{\tau_1,\tau_2\}$ is ``minimal'' in some other sense.  
\end{example}
	
Inspired by \egref{eg:3}, when the convex ordering \eqref{eq:mregular} 
fails, we may define
\begin{equation}\label{eq:m_OBS_setting}
	\begin{aligned}
		U_0(x)\,=\,\poten{\mu_0}(x),\quad U_k(x):=\,\poten{\mu_k}(x)-C_k,
			~~\text{for}~k=1,\cdots,n\,,\\[6pt]
		\text{where}~~C_k:=\,\sup\nolimits_{x\in\R}
			\big\{\poten{\mu_k}(x)-U_{k-1}(x)\big\}.
	\end{aligned}
\end{equation}
Same as before, we assume that the diffusion coefficient $\sigma$ is of 
linear growth and Lipschitz continuous:
\begin{equation}\label{eq:m_Lip_LG}
	\begin{aligned}
		\text{there exists}~L>0,~~\text{s.t.}~~\forall\,x,y\in\R,&\\[6pt]
		|\sigma(x)-\sigma(y)|<&\,L|x-y|,~~|\sigma(x)|<L(1+|x|)\,;
	\end{aligned}
\end{equation}
Moreover, we impose the assumption, similar as \eqref{eq:ellipticity},
\begin{equation}\label{eq:m_ellipticity}
	\begin{aligned}
		\text{for each compact}~K\subset\big\{x:\,\exists\, 
			k=1,\cdots,n,~\text{s.t.}~U_{k-1}(x)>U_k(x)\big\},\\[6pt]
		\text{there exists}~~C_K>0,~~\text{s.t.}~~\forall\,x\in K,
			~~\sigma(x)\,\geq\,C_K\,>\,0\,.
	\end{aligned}
\end{equation}
	
Consider the iterated obstacle problems \mOBS$(\sigma,U_0,u_{k-1}+U_k
-U_{k-1})$ as follows, for $k=1,\cdots,n$,
\begin{equation}\tag{$\mathrm{OBS}$}\label{eq:mOBS}
	\min\Big\{\,\mathrm{L}u_k,\ \big(u_k-u_{k-1}\big)
		-\big(U_k-U_{k-1}\big)\Big\}=0,
		\quad u_k(0,\cdot)=U_0(\cdot)
\end{equation}
where $u_0(t,\cdot):=U_0(\cdot)$ for all $t\geq0$, and the operator 
$\mathrm{L}$ is defined as before.
	
Given the viscosity solutions $u_k$ to \mOBS$(\sigma,U_0,u_{k-1}+U_k
-U_{k-1})$ for $k=1,\cdots,n$, one can define
\begin{equation}\label{eq:defn_tau_k}
	\begin{aligned}
		\tau_0:=0,\quad&\tau_k:=\inf\big\{t\geq\tau_{k-1}:\,
			(t,X_t)\notin D_{k}\big\},
			\quad\text{for}~~k=1,2,\cdots,n\,,\\[6pt]
		\text{where}\quad&D_{k}:=\big\{(t,x):\,u_{k}(t,x)
			-u_{k-1}(t,x)>U_k(x)-U_{k-1}(x)\big\}\,.
	\end{aligned}
\end{equation}	
	
As in \secref{sec:construction}, since $u_{k-1}(0,\cdot)=U_0$, we firstly 
note that the viscosity solution $u_k$ to this \mOBS\ is also the value 
function of the following optimal stopping problem (c.f. 
\citealp[Sect.\,3.4.9]{BensoussanLions:1982}):
\begin{align*}
	\textit{Given that}\quad \dif Y_t\ =&\ \sigma(Y_t)\dif W_t,\quad
		u_k(t,x)\,=\,\sup\nolimits_{\theta\leq t}J^k_{t,x}(\theta),\\[6pt]
	\textit{\ \ \,where}\quad J^k_{t,x}(\theta)\,:=&\ \E^x\big[U_0(Y_\theta)
		\ind_{\theta=t}+\big(u_{k-1}(t-\theta,\cdot)
		+U_k-U_{k-1}\big)(Y_\theta)\ind_{\theta<t}\big]\\[6pt]
	=&\ \E^x\big[u_{k-1}(t-\theta,Y_\theta)+
		\big(U_k(Y_\theta)-U_{k-1}(Y_\theta)\big)\ind_{\theta<t}\big].
\end{align*}
	
In this section, we will generalize \thmref{thm:mainResult} as follows. 
	
\begin{thm}\label{thm:m_mainResult}
	Suppose that \eqref{eq:m_OBS_setting}-\eqref{eq:m_ellipticity} hold, 
	let the sequence of stopping times $\{\tau_k\}$ be given by
	\eqref{eq:defn_tau_k}. Then, for all $k=1,\cdots,n$, we have that
	\begin{equation}\label{eq:each_k}
		X_{\tau_k}\sim\mu_{k},\quad u_k(t,x)=-\E^{\mu_0}\big|x
			-X_{t\wedge\tau_k}\big|,\quad u_k(t,\cdot)\searrow 
			U_k~\text{ as }t\to\infty.
	\end{equation}
\end{thm}
	
Obviously, the desired result directly follows from \thmref{thm:mainResult} 
for $k=1$. Next we will prove \thmref{thm:m_mainResult} by induction.
Firstly, using this connection between obstacle problems and optimal 
stopping problems, we will generalize \lemref{lem:prepare} and 
\lemref{lem:u_distribution} as follows.

\begin{lem}\label{lem:k_barrier}
	Suppose that \eqref{eq:m_OBS_setting}-\eqref{eq:m_ellipticity} hold, and 
	moreover, \eqref{eq:each_k} holds for some $k=1,\cdots,n-1$, then 
	\begin{enumerate}[i), leftmargin = 2.1em, labelsep = .9em ]
		
		\vspace*{-2pt}
		\item\label{itm:k_decrease} $u_{k+1}$, $u_{k+1}-u_{k}$ are 
			non-increasing in $t$, and $D_{k+1}$ is a Root's barrier;
		
		\vspace*{-2pt}
		\item\label{itm:k_distribution} there exist probability 
			distributions $\mu_{k+1}^t$ such that $u_{k+1}(t,\cdot)
			=\poten{\mu^t_{k+1}}$.
			
	\end{enumerate}
\end{lem}	
	
\begin{proof}
	For $s\leq t$, define $\widetilde{u}_k^t(r,x):=u_k(t-r,x)$ and 
	$\widetilde{w}_k^{t,s}(t,x):=(\widetilde{u}_k^t-
	\widetilde{u}_k^s)(r,x)$. Given the assumption that \eqref{eq:each_k} 
	holds for $k$, it follows from \citet[Lem.\,5.2]{CoxOblojTouzi:2018} 
	that both $\big\{\widetilde{u}_{k}^t(r,Y_r)\big\}_{r\leq t}$ and 
	$\big\{\widetilde{w}_{k}^{t,s}(r,Y_r)\big\}_{r\leq s}$ are 
	supermartingales where the process $Y$ is an independent copy of $X$.
		
	For any stopping time $\theta\leq t$, we then have that
	\begin{align*}
		J_{t,x}^{k+1}(\theta)\,&-\,J_{s,x}^{k+1}(s\wedge\theta)\\[6pt] 
		&=\ \E^x\big[ \widetilde{u}_{k}^t(\theta,Y_\theta)
			-\widetilde{u}_{k}^s(s\wedge\theta,Y_{s\wedge\theta})\big]
			+\E^x\big[\big(U_{k+1}(Y_{\theta})-U_k(Y_{\theta})\big)
			\ind_{s\leq\theta<t}\big]\\[6pt]
		\ &\leq\ \E^x\big[\widetilde{w}_{k}^{t,s}(s\wedge\theta,
			Y_{s\wedge\theta})\big]\ \leq\ \widetilde{w}_{k}^{t,s}(0,x)
			\ =\ u_{k}(t,x)-u_{k}(s,x),
	\end{align*}
	where the inequalities hold because $\{\widetilde{u}_{k}^t(r,X_r)\}$ 
	and $\{\widetilde{w}_{k}^{t,s}(r,X_r)\}$ are supermartingales, and 
	$U_{k+1}\leq U_{k}$. Taking supremum over $\theta\leq t$, we conclude 
	the non-increase of $u_{k+1}-u_{k}$ in $t$:
	\begin{equation}\label{eq:w_decreasing}
		\begin{aligned}
			u_{k+1}(t,x)\,-\,u_{k}(t,x)\,&=\,\sup\nolimits_{\theta\leq t}
				J_{t,x}^{k+1}(\theta)-u_{k}(t,x)\\[6pt]
			&\leq\,\sup\nolimits_{\theta\leq t}J_{s,x}^{k+1}(s\wedge\theta)
				-u_{k}(s,x)\,\leq\,u_{k+1}(s,x)-u_{k}(s,x).
		\end{aligned}
	\end{equation}
	Therefore $D_{k+1}$ is a Root's barrier, and $u_{k+1}$ inherits from 
	$u_{k}$ the non-increase in $t$. At last, \itmref{itm:k_distribution} 
	follows from a similar proof as \lemref{lem:u_distribution}.
\end{proof}

Now we can prove the main result of this section.
	
\begin{proof}[Proof of \thmref{thm:m_mainResult}]
	As mentioned above, the desired results hold for $k=1$. For general 
	$k=1,\cdots,n-1$, suppose that \eqref{eq:each_k} holds for $k$.
		
	Since $u_{k+1}-u_k\geq U_{k+1}-U_k$ and $u_{k}\geq U_{k}$, we have that 
	$u_{k+1}\geq U_{k+1}$. Then we can define $\widehat{U}_{k+1}(x):=
	\lim_{t\to\infty}u_{k+1}(t,x)\geq U_{k+1}(x)$, and hence, there exists 
	some constant $C_L$ and a measure $\widehat{\mu}_{k+1}$ such that 
	$\mu_{k+1}^t\Longrightarrow\widehat{\mu}_{k+1}$ and $\widehat{U}_{k+1}
	=\poten{\widehat{\mu}_{k+1}}-C_L$ on $\R$. We also define
	\begin{align*}
		\widehat{D}_{k+1}\,:=\,\big\{(t,x):u_{k+1}(t,x)-u_k(t,x)
			>\widehat{U}_{k+1}(x)-U_k(x)\big\}\,&\subset\,D_{k+1},\\[6pt]
		\text{and}\quad\ \widehat{\tau}_{k+1}
			\,:=\,\inf\big\{t>\tau_k:(t,X_t)\notin\widehat{D}_{k+1}\big\}
			\,&\leq\,\tau_{k+1}.
	\end{align*}
		
	Fix some $t>0$. Let $s=0$ in \eqref{eq:w_decreasing}, we have that
	$u_{k+1}(t,x)\leq u_{k}(t,x)$. Moreover, according to
	\lemref{lem:k_barrier}, we then have that $\mu_k^t=\mathcal{L}
	(X_{t\wedge\tau_k})$ and $\mu_{k}^t\preceq\mu_{k+1}^t$ in convex order.
	Define $v_j(s,x):=u_j(s\wedge t,x)$ for $j=k,k+1$. One can verify that
	\begin{align*}
		\min\big\{(\mathrm{L}v_{k+1})(s,x),~~(v_{k+1}-v_{k})(s,x)
			-(u_{k+1}-u_{k})(t,x)\big\}\,=\,0,
	\end{align*}
	which implies that 
	\begin{equation*}
		u_{k+1}(s\wedge t,x) = \sup\nolimits_{\theta\leq s}\E^x\big[
			u_{k}((s-\theta)\wedge t,Y_\theta)+\big(u_{k+1}(t,Y_\theta)
			-u_{k}(t,Y_\theta)\big)\ind_{\theta<s}\big].
	\end{equation*}
	Define $\tau_{k+1}^t:=\inf\big\{s\geq t\wedge\tau_{k}:(s,X_s)\notin
		D_{k+1}^t\big\}$, where
	\begin{equation*}
		D_{k+1}^t\,:=\,\Big\{(s,x):u_{k+1}(s\wedge t,x)
			-u_{k}(s\wedge t,x)>u_{k+1}(t,x)-u_{k}(t,x)\Big\},
	\end{equation*}
	According to \citet[Thm.\,4.1]{CoxOblojTouzi:2018}, we have that
	$\tau_{k+1}^t\geq t\wedge\tau_k$, $X_{\tau_{k+1}^t}\sim \mu_{k+1}^t$.
	Similar as \eqref{eq:tau_limit}, we can show that, as $t\to\infty$, \,
	$D^t_{k+1}\nearrow \widehat{D}_{k+1}$, and $\tau^t_{k+1}\nearrow
	\widehat{\tau}_{k+1}$, and hence, $X_{\widehat{\tau}_{k+1}}
	=\lim_{t\to\infty}X_{\tau_{k+1}^t}\sim\widehat{\mu}_{k+1}$.
	Then, by a simple modification of the proof of \thmref{thm:mainResult}, 
	we can conclude that \eqref{eq:each_k} also holds for $k+1$.
	This completes our proof.
\end{proof}

\subsection{Minimality of Root's embeddings to multi-marginal SEP}
\label{sec:multi_minimality}
	
Let $\mtau=\{\tau_1,\cdots\tau_n\}$ be the sequence given by 
\eqref{eq:defn_tau_k}. When the convex ordering condition 
\eqref{eq:mregular} holds, we have that $C_k=0$ for all $k$ in
\eqref{eq:m_OBS_setting}. According to 
\citet[Thm.\,3.1]{CoxOblojTouzi:2018}, $\tau_k$ are UI stopping times. 
However, as mentioned before, we cannot expect so in the absence of 
\eqref{eq:mregular}. Same as \secref{sec:minimality}, we now consider the 
minimality of our solution to \mSEP$(\sigma,\mu_0,\mmu)$.
	
The first hitting time $\tau_1$ is a minimal Root's stopping time under 
$\pb^{\mu_0}$ (\thmref{thm:minimalResult}).	However, as seen in 
\egref{eg:3}, the subsequent stopping times do not inherit the minimality
(unless $\mu_0\preceq\cdots\preceq\mu_{k-1}$ in convex order, or 
equivalently, $\tau_1,\cdots,\tau_{k-1}$ are UI stopping times). Hence, 
we focus on the ``minimality'' of a sequence of stopping times in some 
other sense.
	
\begin{defn}[Minimal sequence of stopping times]
\label{defn:minimal_sequence}
	A non-decreasing sequence of stopping times $\mtau=\{\tau_k\}_{k=1,
	\cdots,n}$ for the process $X$ is minimal if whenever 
	$\mtheta=\{\theta_k\}_{k=1,\cdots,n}$ is a non-decreasing sequence of 
	stopping times such that $\theta_{k}\leq\tau_{k}$ and
	$\mathcal{L}(X_{\theta_{k}})=\mathcal{L}(X_{\tau_{k}})$ for all $k$ 
	then $\mtau=\mtheta$ almost surely.
	
	Moreover we say that $\mtau$ is a minimal solution to \mSEP$(\sigma,
	\mu_0,\mmu)$ if $\mtau$ is a minimal sequence and a solution to
	\mSEP$(\sigma,\mu_0,\mmu)$ simultaneously.
\end{defn}
	
\begin{prop}\label{prop:minimal_sequence_stopping}
	Denote $\tau_0=0$. A non-decreasing sequence $\mtau$ is a minimal 
	sequence of stopping times if and only if 
	\begin{equation}\label{eq:stopping2sequence}
		\begin{aligned}
			\text{ whenever $\,\theta$ is a stopping time such that }
				\,\tau_{k-1}\leq\theta\leq\tau_k&\\[4pt]
			\text{ and }\mathcal{L}(X_{\theta})
				=\mathcal{L}(X_{\tau_{k}})\ \text{for some $k$ then }
				\theta=\tau_k&\,.
		\end{aligned}
	\end{equation}
\end{prop}

\begin{proof}
	Firstly suppose that the sequence $\mtau$ satisfies 
	\eqref{eq:stopping2sequence}, and $\mtheta$ is a sequence such that 
	$\theta_{k}\leq\tau_k$ and $\mathcal{L}(X_{\theta_{k}})=
	\mathcal{L}(X_{\tau_{k}})$. Obviously $\theta_1=\tau_1$ since
	$\tau_1$ is a minimal stopping time in the sense of Monroe. It then 
	follows from \eqref{eq:stopping2sequence} and induction that 
	$\theta_k=\tau_k$ for all $k$.
		
	Conversely, suppose that $\mtau$ is a minimal sequence. For any 
	$k^\ast\in\{1,\cdots n\}$, let $\theta$ be a stopping time such that 
	$\tau_{k^\ast-1}\leq\theta\leq\tau_{k^\ast}$ and $\mathcal{L}
	(X_{\theta})=\mathcal{L}(X_{\tau_{k^\ast}})$. Replacing $\tau_{k^\ast}$ 
	by $\theta$ in $\mtau$, one then have another sequence which is
	also non-decreasing and embeds same marginal distributions as 
	$\{\tau_k\}$ does. It then follows that $\theta=\tau_{k^\ast}$ since 
	$\mtau$ is a minimal sequence.
\end{proof}

Now we focus on the property described in \eqref{eq:stopping2sequence}.
Given a pair of stopping times $S\leq T$, we say that $T$ is \emph{minimal 
with respect to $S$} if whenever $R$ is a stopping time s.t. $S\leq R\leq T$ 
and $\mathcal{L}(X_R)=\mathcal{L}(X_T)$ then $R=T$ a.s. (as described in
\eqref{eq:stopping2sequence}). By a similar proof as in 
\citet[Prop.\,2]{Monroe:1972a}, for any stopping time $R\geq S$ there is a 
stopping time $T\leq R$ which is minimal with respect to $S$ and embeds 
same distribution as $R$ does. Further, by a careful review and simple 
modification of the arguments in \citet[Sect.\,2]{CoxHobson:2006}, 
\citet[Appx.\,A,\,B]{Cox:2008} and \secref{sec:minimality} of this work,
one can generalize \thmref{thm:CoxMinimality} and \thmref{thm:minimal} as 
follows:
	
\begin{prop}\label{prop:m_Cox_minimal}
	Suppose that $S\leq T$ are stopping times such that $X_{S}\sim\nu$, 
	$X_T\sim\mu$ under some probability measure $\pb$. The set 
	$\mathcal{A}$ and its upper/lower bound $a_\pm$ are defined as in 
	\eqref{eq:defn_A&C} and \eqref{eq:defn_a_pm}. Moreover, for the
	set $\mathcal{A}$ and some horizontal level $\gamma$, denote the first 
	hitting times after $S$ by $H^{S}_\mathcal{A}$ and $H^{S}_\gamma$:
	\begin{equation*}
		H^{S}_\mathcal{A}\,=\,\inf\big\{t\geq S:X_t\in\mathcal{A}\big\}
			\quad\text{and}\quad H^{S}_\gamma\,=\,\inf\big\{t\geq S:
			X_t=\gamma\big\}.
	\end{equation*}
	Then the following statements are equivalent:
	
	\vspace{-3pt}
	\begin{enumerate}[i), leftmargin = 2.4em, labelsep = .9em ]
		
		\item\label{itm:m_minimal_1} $T$ is minimal with respect to $S$;
		
		\item\label{itm:m_minimal_3} $T\leq H^{S}_{\mathcal{A}}$ and for 
			all stopping times $R$ such that $S\leq R\leq T$, 
			\begin{equation*}
				\E\big[X_T\big|\mathcal{F}_R\big]
					\leq X_{R}\,\text{ on }\,\{X_{S}\geq a_-\};\qquad
				\E\big[X_T\big|\mathcal{F}_R\big]
					\geq X_{R}\,\text{ on }\,\{X_{S}\leq a_+\};
			\end{equation*}
		
		\item\label{itm:m_minimal_5} $T\leq H^{S}_{\mathcal{A}}$ and as 
			$\gamma\rightarrow\infty$, 
			\begin{equation*}
				\gamma\,\pb\big[T>H^{S}_{-\gamma},\ X_{S}\geq a_-\big]
					\longrightarrow 0;\qquad
				\gamma\,\pb\big[T>H^{S}_{+\gamma},\ X_{S}\leq a_+\big]
					\longrightarrow 0;
			\end{equation*}
		
		\item \label{itm:m_minimal_6} $\inf_{x\in \R}\E\big[L_T^x
			-L_{S}^x\big]\,=\,0$.
	\end{enumerate}
	Further, if $\,\exists\,a\in\R$ s.t. $\pb[T\leq H^{S}_a]=1$, $T$ is 
	minimal with respect to $S$.

\end{prop}
	
Given a solution to \mSEP$(\sigma,\mu_0,\mmu)$, denoted by $\mtheta=
\{\theta_k\}$, same as in the proof of \lemref{lem:MinStop_2_MinC} we have 
that there exists $\{c_k\}$ such that $u^\mtheta_k=-\E^{\mu_0}\left[x
-X_{t\wedge\theta_k}\right]\to\poten{\mu_k}(x)-c_k$ and
\begin{equation*}
	c_k\,\equiv\,\E^{\mu_0}\left[\int_0^{\theta_k}\sgn(x-X_s)
		\dif X_s\right],\quad\text{for all}~x\in\R.
\end{equation*}
And then it follows from the definition of $\{C_k\}$ (recall 
\eqref{eq:m_OBS_setting})
\begin{align*}
	C_k-C_{k-1}\,&=\,\sup\nolimits_{x\in\R}\big\{\poten{\mu_k}(x)
		-\poten{\mu_{k-1}}(x)\big\}\\[6pt]
	\,&=\,\big(c_{k}-c_{k-1}\big)-\inf\nolimits_{x\in\R}
		\E^{\mu_0}\big[L^x_{\theta_k}-L^x_{\theta_{k-1}}\big].
\end{align*}
Noting that $C_{0}=c_0=0$, by \propref{prop:minimal_sequence_stopping}
\& \ref{prop:m_Cox_minimal} we then have that
\begin{align*}
	C_k\,=\,c_k,\quad\text{for all}~k\quad&\Longleftrightarrow\quad
		\inf\nolimits_{x\in\R}\E^{\mu_0}\big[L^x_{\theta_k}-
		L^x_{\theta_{k-1}}\big]\,=\,0,\quad\text{for all}~k\\[6pt]
	&\Longleftrightarrow\quad\,\mtheta
		~\text{is a minimal sequence of stopping times.}
\end{align*}
This result extends \thmref{thm:minimalResult} to the multi-marginal \mSEP\
as follows.
	
\begin{thm}\label{thm:m_minimalResult}
	Suppose that \eqref{eq:m_OBS_setting}-\eqref{eq:m_ellipticity} hold,
	then 
	\begin{enumerate}[i), leftmargin = 2.1em, labelsep = .9em ]

		\item\label{itm:minimal_sequence} The sequence $\mtau$ given by 
			\eqref{eq:defn_tau_k} is a minimal solution
			to \mSEP$(\sigma,\mu_0,\mmu)$.
		
		\item\label{itm:m_minimal_condition} Let $\mtheta$ be a minimal 
			solution to \mSEP$(\sigma,\mu_0,\mmu)$, then the potential 
			process $u^\mtheta_k(t,x)\searrow U_k(x)$ as $t\to\infty$ 
			for all $k=1,\cdots,n$.
		
	\end{enumerate}
\end{thm}

\subsection{Optimality of Root's embeddings to multi-marginal SEP}
\label{sec:multi_optimality}
	
At last, we consider the optimality obtained in \secref{sec:optimality} to
the multi-marginal distributions cases. Suppose that $\mtau=\{\tau_k\}$
and $\mtheta=\{\theta_k\}$ are \emph{minimal} solutions to \mSEP$(\sigma,
\mu_0,\mmu)$, among which $\mtau$ is given by \eqref{eq:defn_tau_k}. As 
before, define the potential processes
\begin{equation*}
	u^\mtheta_k(t,x)\,=\,-\E^{\mu_0}\big|x-X_{t\wedge\theta_{k} }\big|,
	\qquad
	u^\mtau_k(t,x)\,=\,-\E^{\mu_0}\big|x-X_{t\wedge\tau_{k} }\big|.
\end{equation*}
According to \thmref{thm:m_minimalResult}, $u^\mtheta_j(t,x)\searrow 
\poten{\mu_j}(x)-C_j=U_j(x)$ as $t\to\infty$ for $j=k-1,k$, where $C_j$ 
and $U_j$ are as defined in \eqref{eq:m_OBS_setting}. Moreover, one can easily verify that, for $t\geq s\geq 0$,
\begin{align*}
	\big[u^\mtheta_k(t,x)-u^\mtheta_{k-1}(t,x)\big]
		&-\big[u^\mtheta_k(s,x)-u^\mtheta_{k-1}(s,x)\big]\\[6pt]
	\,&=\,-\E^{\mu_0}\big[ \big(L^x_{t\wedge\theta_{k}}
		-L^x_{t\wedge\theta_{k-1}}\big)-\big(L^x_{s\wedge\theta_{k}}
		-L^x_{s\wedge\theta_{k-1}}\big)\big]\,\leq\,0.
\end{align*}
It follows that $u^\mtheta_k(t,x)-u^\mtheta_{k-1}(t,x)\searrow 
U_k(x)-U_{k-1}(x)$ as $t\to\infty$, and hence, 
\begin{align*}
	\min\Big\{\,(\mathrm{L}u^\mtheta_k)(t,x)&,\ 
		\big(u^\mtheta_k-u^\mtheta_{k-1}\big)(t,x)
		-\big(U_k-U_{k-1}\big)(x)\Big\}\\[6pt]
	&\geq\,0\, =\, \min\Big\{\,(\mathrm{L}u^\mtau_k)(t,x),
		\ \big(u^\mtau_k-u^\mtau_{k-1}\big)(t,x)
		-\big(U_k-U_{k-1}\big)(x)\Big\}
\end{align*}
where $\mathrm{L}u^\mtheta_k\geq 0$ follows from 
\citeauthor{GassiatOberhauserReis:2015}(see \eqref{eq:viscosity_super} 
\& the footnote on \ppref{eq:viscosity_super}).
	
Suppose that $u^\mtau_{k-1}(t,x)\leq u^\mtheta_{k-1}(t,x)$ for all $(t,x)$. 
Since $u^\mtheta_k$, $u^\mtau_{k}$, $u^\mtau_{k-1}$, $U_k$, $U_{k-1}$ are 
all Lipschitz continuous in $x$ (uniformly in $t$), a slight extension of 
\citet[Thm.\,5]{GassiatOberhauserReis:2015} implies that  $u^\mtau_k\leq u^
\mtheta_{k}$. Since $u^\mtau_1\leq u^\mtheta_1$ (\propref{prop:optimality}), 
we have the following result by induction.
	
\begin{prop}\label{prop:k_optimality}
	Suppose that \eqref{eq:m_OBS_setting}-\eqref{eq:m_ellipticity} hold.
	Let $\mtau$ and $\mtheta$ be minimal solutions to \mSEP$(\sigma,
	\mu_0,\mmu)$, among which $\mtau$ is given by \eqref{eq:defn_tau_k}. 
	Then $\,u^{\mtau}_k(t,x) \leq u^\mtheta_k(t,x)$ for all 
	$(k,t,x)\in\{1,\cdots,n\}\times[0,+\infty)\times\R.$
\end{prop}

Now we follow the work of \citet{CoxOblojTouzi:2018} and 
\secref{sec:optimality} of this work to get the multi-marginal analogue of
\thmref{thm:path_optimality}. 
	
\begin{thm}\label{thm:m_path_optimality}
	Suppose that \eqref{eq:sigma_condition} and 
	\eqref{eq:m_OBS_setting}-\eqref{eq:m_ellipticity} hold. Let $\mtau$ 
	and $\mtheta$ be the minimal solutions to \mSEP$(\sigma,\mu_0,\mmu)$,
	among which $\mtau$ is given by \eqref{eq:defn_tau_k}. Then $\E^{\mu_0}
	\big[F(\tau_n)\big]\geq\E^{\mu_0}\big[F(\theta_n)\big]$ for all 
	non-decreasing concave function $F$.
\end{thm}
	
\begin{proof}
	Without loss of generality, we always assume that $F(0)=0$. In addition, 
	we may firstly suppose that \eqref{eq:vanish_f} holds for the right 
	derivative $f=F'_+$, and the general case follows from monotone 
	convergence theorem (see the last paragraph of the proof of 
	\thmref{thm:path_optimality}). 
		
	Let $\{D_k\}$ be the sequence of barriers given by 
	\eqref{eq:defn_tau_k}. For $(k,t,x)\in\{1,\cdots,n\}\times[0,+\infty)
	\times\R$, we define stopping times $\eta_k=\inf\{s\geq t:\,(s,X_s)
	\notin D_k\}$ under $\pb^{(t,x)}$. Further, define $M_{n+1}(t,x):=f(t)$ 
	and
	\begin{equation*}
		M_{k}(t,x):=\E^{(t,x)}\big[M_{k+1}(\eta_k,X_{\eta_k})\big],
			\quad\text{for}~\,k=1,\cdots,n.
	\end{equation*}
	Similar as in the proof of \citet[Lem.\,3.4]{CoxOblojTouzi:2018},
	$M_k(t,x)=\E^{(t,x)}\big[f(\zeta^k)\big]$ where $\zeta^k$ is the first
	time we exit $D_n$, having previously exited $D_k,\cdots,D_{n-1}$ in 
	sequence. Hence $\zeta^k$ is non-increasing in $k$, and $M_k$ is 
	non-decreasing in $k$ because $f$ is non-increasing. Moreover, for any 
	$k$, if $(t,x)\notin D_k$, then $\zeta^k=\zeta^{k+1}$, 
	$\pb^{(t,x)}$-a.s.. As conclusion,
	\begin{equation*}
		M_k(t,x)\leq M_{k+1}(t,x),\quad\text{for all}~
			(t,x)\in\R_+\times\R,\quad\text{with ``$=$'' if}~t\geq R_k(x),
	\end{equation*}
	where $R_k$ denotes the boundary function of $D_k$. Hence, given 
	$(t_k,x_k)_{k=0,\cdots,n}$ with $0=t_0\leq t_1\leq\cdots\leq t_n$, we 
	deduce that
	\begin{align*}
		\int_0^{t_k}\!M_k(s,x_k)\dif s\,
			+&\int_0^{R_k(x_k)}\!\big(M_{k+1}-M_{k}\big)(s,x_k)\dif s
			-\!\int_0^{t_{k-1}}\!M_k(s,x_{k-1})\dif s\\[6pt]
		=\ &\int_0^{R_k(x_k)}\!M_{k+1}(s,x_k)\dif s-\!\int_{t_k}^{R_k(x_k)}
			\!M_k(s,x_k)\dif s-\!\int_0^{t_{k-1}}
			\!M_k(s,x_{k-1})\dif s\\[6pt]
		\geq\ &\int_0^{t_k}\!M_{k+1}(s,x_k)\dif s-\!\int_0^{t_{k-1}}
			\!M_k(s,x_{k-1})\dif s,\qquad
			\text{with ``$=$'' if }~t_k\geq R_k(x_k).
	\end{align*}
	Taking sum for $k=1,\cdots,n$ and noting  $F(0)=0$ and $t_0=0$, we 
	have that
	\begin{align*}
		\sum_{k=1}^{n}\Bigg\{&\int_0^{t_k}\!M_k(s,x_k)\dif s
			+\int_0^{R_k(x_k)}\!\big(M_{k+1}-M_{k}\big)(s,x_k)\dif s
			-\int_0^{t_{k-1}}\!M_k(s,x_{k-1})\dif s\Bigg\}\\[6pt]
		\geq\ &\int_0^{t_n}\!M_{n+1}(s,x_k)\dif s
				-\!\int_0^{t_{0}}\!M_1(s,x_{0})\dif s
				\qquad\big(\text{with ``$=$'' if}~t_k\geq R_k(x_k)
				~\text{for all}~k\big)\\[6pt]
		=\ &F(t_n)
	\end{align*}
	Since $\mathcal{L}(X_{\tau_k})=\mathcal{L}(X_{\theta_k})=\mu_k$ 
	and $\tau_k\geq R(X_{\tau_k})$, we then have that
	\begin{align*}
		\E^{\mu_0}\big[F(\tau_{n})\big]
			\ -\ \E^{\mu_0}\big[F(\theta_{n})\big] 
			\,\geq\ \sum_{k=1}^{n}\Bigg\{&\E^{\mu_0}\left[\int_0^{\tau_k}
			 	\!M_k(s,X_{\tau_k})\dif s-\!\int_0^{\tau_{k-1}}
			 	\!M_k(s,X_{\tau_{k-1}})\dif s\right]\\[6pt]
			-\,&\E^{\mu_0}\left[\int_0^{\theta_k}\!M_k(s,X_{\theta_k})
				\dif s-\!\int_0^{\theta_{k-1}}\!M_k(s,X_{\theta_{k-1}})
				\dif s\right]\Bigg\}.
	\end{align*}
	Thus, to see $\E^{\mu_0}\big[F(\tau_{n})-F(\theta_{n})\big]\geq 0$,
	it is sufficient to show that, for all $k$,
	\begin{equation}\label{eq:M_k_ineq}
		\begin{aligned}
			\E^{\mu_0}\bigg[\int_0^{\tau_k}\!M_k(s,X_{\tau_k})\dif s
				\,&-\int_0^{\tau_{k-1}}\!M_k(s,X_{\tau_{k-1}})\dif s
				\bigg]\\[6pt]
			&\geq\ \E^{\mu_0}\bigg[\int_0^{\theta_k}\!M_k(s,X_{\theta_k})
				\dif s-\!\int_0^{\theta_{k-1}}\!M_k(s,X_{\theta_{k-1}})
				\dif s\bigg].
		\end{aligned}
	\end{equation}
		
	We are going to show \eqref{eq:M_k_ineq} in the same manner as proving 
	\eqref{eq:M_compare}. Define
	\begin{align*}
		G_k(t,x)=\!\int_0^t\!M_k(s,x)\dif s-Z_k(x),\quad\text{where\,}~
			Z_k(x)=\int_0^x\!\dif y\!\int_0^y\dfrac{2M_k(0,z)}{\sigma^2(z)}
			\dif z.
	\end{align*}
		
	Same as in the proof of \thmref{thm:path_optimality}, one can show that
	$Z_k(X_t)$ and $G_{k}(t,X_t)$ are integrable according to
	\eqref{eq:sigma_condition} and \eqref{eq:vanish_f}. Then a simple 
	modification of the proof of \citet[Lem.\,A.1]{CoxOblojTouzi:2018} says 
	that $\big\{G_k(t,X_{t})\big\}$ is a $\pb^{\mu_0}$-supermartingale
	and a $\pb^{\mu_0}$-martingale on $\big[\tau_{k-1},\tau_k\big]$.
	Thus, for the sequences $\mtau$ and $\mtheta$, we deduce that (recall 
	\eqref{eq:Z_expectation})\footnote{
		We have mentioned in \rmkref{rmk:unsufficient_comparison_1} that
		the comparison $u^\tau_k\leq u^\theta_k$ alone is not sufficient 
		to show \eqref{eq:M_k_ineq}. In more detail, one may find that we 
		need a stronger result than \propref{prop:k_optimality},
		that is $u^\mtau_k-u^\mtau_{k-1}\leq u^\mtheta_k-u^\mtheta_{k-1}$ 
		everywhere.}
	\begin{align}\label{eq:M_Z_ineq}
		\begin{aligned}
			\E^{\mu_0}\bigg[\int_0^{t\wedge\tau_k}
				&\!M_k(s,X_{t\wedge\tau_k})\dif s-\int_0^{t\wedge\tau_{k-1}}
				\!M_k(s,X_{t\wedge\tau_{k-1}})\dif s\bigg]\\[6pt]
			&\begin{aligned}
				-\E^{\mu_0}\bigg[\int_0^{t\wedge\theta_k}
					\!M_k(s,&X_{t\wedge\theta_k})\dif s-
					\int_0^{t\wedge\theta_{k-1}}
					\!M_k(s,X_{t\wedge\theta_{k-1}})\dif s\bigg]\\[6pt]
				\geq\ &\int_{\R}\dfrac{\big(u^\mtheta_k-
					u^\mtheta_{k-1}\big)(t,y)-\big(u^\mtau_k-
					u^\mtau_{k-1}\big)(t,y)}{2}\, Z_k''(y) \dif y.
			\end{aligned}
		\end{aligned}
	\end{align}	
	By minimality,
	$u^{\mrho}_k(t,\cdot)-u^{\mrho}_{k-1}(t,\cdot)\searrow U_k-U_{k-1}$
	for $\mrho =\mtau,\mtheta$, and then,
	\begin{align*}
		\big(U_k-U_{k-1}\big)-\big(u^\mtau_k-u^\mtau_{k-1}\big)
			\,&\leq\,\big(u^\mtheta_k-u^\mtheta_{k-1}\big)
			-\big(u^\mtau_k-u^\mtau_{k-1}\big)\\[6pt]
		\,&\leq\,\big(u^\mtheta_k-u^\mtheta_{k-1}\big)
			-\big(U_k-U_{k-1}\big).
	\end{align*}
	It follows from monotone convergence and squeeze theorem that the RHS 
	of \eqref{eq:M_Z_ineq} vanishes as $t\to\infty$. Then 
	\eqref{eq:M_k_ineq} follows from dominated convergence on 
	\eqref{eq:M_Z_ineq} (similar as in the proof of 
	\thmref{thm:path_optimality}).
\end{proof}		


\section*{Acknowledgments}
\addcontentsline{toc}{section}{Acknowledgments}

I would like to show my gratitude to the anonymous editors and reviewers
of \emph{Stochastic Processes and their Applications} for the suggestion on 
the extension to multi-marginal SEP and also the encouragement to resubmit 
this manuscript. I am also immensely grateful to Alexander M.~G.~Cox for the 
helpful discussion on this project during my visit at University of Bath.

	
\small

\linespread{1.1}	
	
\addcontentsline{toc}{section}{\refname}

\bibliographystyle{my_plainnat} 
\bibliography{minimal_Root_embedding}

\begin{thebibliography}{30}
\providecommand{\natexlab}[1]{#1}
\providecommand{\url}[1]{\texttt{#1}}
\expandafter\ifx\csname urlstyle\endcsname\relax
  \providecommand{\doi}[1]{doi: #1}\else
  \providecommand{\doi}{doi: \begingroup \urlstyle{rm}\Url}\fi

\bibitem[Az{\'e}ma and Yor(1979)Jacques Az{\'e}ma and Marc Yor]{AzemaYor:1979a}
{Az{\'e}ma, Jacques and Yor, Marc}.
\newblock Une solution simple au probl\`eme de {S}korokhod.
\newblock In \emph{S\'eminaire de {P}robabilit\'es {XIII} ({U}niv.
  {S}trasbourg, {S}trasbourg, 1977/78)}, volume 721 of \emph{Lecture Notes in
  Math.}, pages 90--115. Springer-Verlag, Berlin Heidelberg, 1979.
\newblock ISBN 978-3-540-09505-7.
\newblock \doi{10.1007/BFb0070852}.

\bibitem[Beiglb{\"o}ck et~al.(2017)Mathias Beiglb{\"o}ck, Alexander M.~G. Cox,
  and Martin Huesmann]{BeiglbockCoxHuesmann:2017}
{Beiglb{\"o}ck, Mathias and Cox, Alexander M. G. and Huesmann, Martin}.
\newblock Optimal transport and {S}korokhod embedding.
\newblock \emph{Invent. Math.}, 208\penalty0 (2):\penalty0 327--400, 2017.
\newblock ISSN 0020-9910.
\newblock \doi{10.1007/s00222-016-0692-2}.

\bibitem[Bensoussan and Lions(1982)Alain Bensoussan and Jacques-Louis
  Lions]{BensoussanLions:1982}
{Bensoussan, Alain and Lions, Jacques-Louis}.
\newblock \emph{Applications of variational inequalities in stochastic
  control}, volume~12 of \emph{Studies in Mathematics and its Applications}.
\newblock North-Holland Publishing Co., Amsterdam-New York, first edition,
  1982, pages xi+564.
\newblock ISBN 978-0-444-86358-4.
\newblock \doi{10.1016/s0168-2024(08)x7013-4}.
\newblock Translated from the French.

\bibitem[Chacon(1977)Rafael~V. Chacon]{Chacon:1977}
{Chacon, Rafael V.}
\newblock Potential processes.
\newblock \emph{Trans. Amer. Math. Soc.}, 226:\penalty0 39--58, 1977.
\newblock ISSN 0002-9947.
\newblock \doi{10.1090/S0002-9947-1977-0501374-5}.

\bibitem[Chacon and Walsh(1976)Rafael~V. Chacon and John~B.
  Walsh]{ChaconWalsh:1976}
{Chacon, Rafael V. and Walsh, John B.}
\newblock One-dimensional potential embedding.
\newblock In \emph{S\'eminaire de {P}robabilit\'es {X} ({P}r\`emiere partie,
  {U}niv. {S}trasbourg, {S}trasbourg, ann\'ee universitaire 1974/75)}, volume
  511 of \emph{Lecture Notes in Math.}, pages 19--23. Springer-Verlag, Berlin
  Heidelberg, 1976.
\newblock ISBN 978-3-540-07681-0.
\newblock \doi{10.1007/BFb0101093}.

\bibitem[Chacon(1985)Rene~M. Chacon]{ChaconRene:Thesis}
{Chacon, Rene M.}
\newblock \emph{Barrier stopping times and the filling scheme}.
\newblock PhD thesis, Department of Mathematics, University of Washington, Ann
  Arbor, MI, 1985.
\newblock ProQuest document ID:\
  \href{http://search.proquest.com/docview/303415878?accountid=11719}{\ttfamily
  303415878}.

\bibitem[Cox(2008)Alexander M.~G. Cox]{Cox:2008}
{Cox, Alexander M. G.}
\newblock Extending {C}hacon-{W}alsh: minimality and generalised starting
  distributions.
\newblock In \emph{S\'eminaire de probabilit\'es {XLI}}, volume 1934 of
  \emph{Lecture Notes in Math.}, pages 233--264. Springer-Verlag, Berlin
  Heidelberg, 2008.
\newblock ISBN 978-3-540-77912-4.
\newblock \doi{10.1007/978-3-540-77913-1_12}.

\bibitem[Cox and Hobson(2006)Alexander M.~G. Cox and David~G.
  Hobson]{CoxHobson:2006}
{Cox, Alexander M. G. and Hobson, David G.}
\newblock Skorokhod embeddings, minimality and non-centred target
  distributions.
\newblock \emph{Probab. Theory Related Fields}, 135\penalty0 (3):\penalty0
  395--414, 2006.
\newblock ISSN 0178-8051.
\newblock \doi{10.1007/s00440-005-0467-y}.

\bibitem[Cox and Wang(2013)Alexander M.~G. Cox and Jiajie Wang]{CoxWang:2013a}
{Cox, Alexander M. G. and Wang, Jiajie}.
\newblock Root's barrier: construction, optimality and applications to variance
  options.
\newblock \emph{Ann. Appl. Probab.}, 23\penalty0 (3):\penalty0 859--894, 2013.
\newblock ISSN 1050-5164.
\newblock \doi{10.1214/12-AAP857}.

\bibitem[Cox et~al.(2018)Alexander M.~G. Cox, Jan Ob{\l}{\'o}j, and Nizar
  Touzi]{CoxOblojTouzi:2018}
{Cox, Alexander M. G. and Ob{\l}{\'o}j, Jan and Touzi, Nizar}.
\newblock The {R}oot solution to the multi-marginal embedding problem: an
  optimal stopping and time-reversal approach.
\newblock \emph{Probab. Theory Related Fields}, 170, Feb 2018.
\newblock ISSN 0178-8051.
\newblock \doi{10.1007/s00440-018-0833-1}.

\bibitem[Dubins(1968)Lester~E. Dubins]{Dubins:1968}
{Dubins, Lester E.}
\newblock On a theorem of {S}korohod.
\newblock \emph{Ann. Math. Statist.}, 39\penalty0 (6):\penalty0 2094--2097,
  1968.
\newblock ISSN 0003-4851.
\newblock \doi{10.1214/aoms/1177698036}.

\bibitem[Dupire(2005)Bruno Dupire]{Dupire:2005}
{Dupire, Bruno}.
\newblock Arbitrage bounds for volatility derivatives as free boundary problem.
\newblock Presentation at \emph{PDE and Mathematical Finance}, KTH Royal
  Institute of Technology, Stockholm, August 2005.
\newblock \sc url:
  \tt\url{http://www.math.kth.se/pde_finance/presentations/Bruno.pdf}\normalfont.

\bibitem[Ekstr\"om and Tysk(2011)Erik Ekstr\"om and Johan
  Tysk]{EkstromTysk:2011}
{Ekstr\"om, Erik and Tysk, Johan}.
\newblock Boundary behaviour of densities for non-negative diffusions.
\newblock Uppsala University, May 2011.
\newblock \sc url: \tt\url{http://www2.math.uu.se/~ekstrom/pq.pdf}\normalfont.

\bibitem[El~Karoui et~al.(1997)Nicole El~Karoui, Christophe Kapoudjian,
  \'{E}tienne Pardoux, Shi~Ge Peng, and Marie-Claire
  Quenez]{ElKarouiPardouxPeng:1997}
{El Karoui, Nicole and Kapoudjian, Christophe and Pardoux, \'{E}tienne and
  Peng, Shi Ge and Quenez, Marie-Claire}.
\newblock Reflected solutions of backward {SDE}'s, and related obstacle
  problems for {PDE}'s.
\newblock \emph{Ann. Probab.}, 25\penalty0 (2):\penalty0 702--737, 1997.
\newblock ISSN 0091-1798.
\newblock \doi{10.1214/aop/1024404416}.

\bibitem[Gassiat et~al.(2015)Paul Gassiat, Harald Oberhauser, and Gon{\c{c}}alo
  dos Reis]{GassiatOberhauserReis:2015}
{Gassiat, Paul and Oberhauser, Harald and dos Reis, Gon{\c{c}}alo}.
\newblock Root's barrier, viscosity solutions of obstacle problems and
  reflected {FBSDE}s.
\newblock \emph{Stochastic Process. Appl.}, 125\penalty0 (12):\penalty0
  4601--4631, 2015.
\newblock ISSN 0304-4149.
\newblock \doi{10.1016/j.spa.2015.07.010}.

\bibitem[Hobson(2011)David~G. Hobson]{Hobson:2011}
{Hobson, David G.}
\newblock The {S}korokhod embedding problem and model-independent bounds for
  option prices.
\newblock In \emph{Paris-{P}rinceton {L}ectures on {M}athematical {F}inance
  2010}, volume 2003 of \emph{Lecture Notes in Math.}, pages 267--318.
  Springer-Verlag, Berlin Heidelberg, 2011.
\newblock ISBN 978-3-642-14659-6.
\newblock \doi{10.1007/978-3-642-14660-2_4}.

\bibitem[It\^o and McKean~Jr.(1974)Kiyosi It\^o and Henry~P.
  McKean~Jr.]{ItoMcKean:1974}
{It\^o, Kiyosi and McKean Jr., Henry P.}
\newblock \emph{Diffusion processes and their sample paths}.
\newblock Classics in Mathematics. Springer-Verlag, Berlin Heidelberg, 1974,
  pages xv+321.
\newblock ISBN 978-3-540-60629-1.
\newblock \doi{10.1007/978-3-642-62025-6}.
\newblock Second printing, corrected, Die Grundlehren der mathematischen
  Wissenschaften, Band 125.

\bibitem[Kiefer(1972)James~E. Kiefer]{Kiefer:1972}
{Kiefer, James E.}
\newblock Skorohod embedding of multivariate {RV}'s, and the sample {DF}.
\newblock \emph{Z. Wahrscheinlichkeitstheorie und Verw. Gebiete}, 24\penalty0
  (1):\penalty0 1--35, 1972.
\newblock ISSN 0044-3719.
\newblock \doi{10.1007/BF00532460}.

\bibitem[Lange(2010)Kenneth~L. Lange]{Lange:2010}
{Lange, Kenneth L.}
\newblock \emph{Applied probability}.
\newblock Springer Texts in Statistics. Springer-Verlag, New York, NY, second
  edition, 2010, pages xvi+436.
\newblock ISBN 978-1-4419-7164-7.
\newblock \doi{10.1007/978-1-4419-7165-4}.

\bibitem[Loynes(1970)Robert~M. Loynes]{Loynes:1970}
{Loynes, Robert M.}
\newblock Stopping times on {B}rownian motion: {S}ome properties of {R}oot's
  construction.
\newblock \emph{Z. Wahrscheinlichkeitstheorie und Verw. Gebiete}, 16\penalty0
  (3):\penalty0 211--218, 1970.
\newblock ISSN 0044-3719.
\newblock \doi{10.1007/BF00534597}.

\bibitem[Monroe(1972)Itrel~E. Monroe]{Monroe:1972a}
{Monroe, Itrel E.}
\newblock On embedding right continuous martingales in {B}rownian motion.
\newblock \emph{Ann. Math. Statist.}, 43\penalty0 (4):\penalty0 1293--1311,
  1972.
\newblock ISSN 0003-4851.
\newblock \doi{10.1214/aoms/1177692480}.

\bibitem[Ob{\l}{\'o}j(2004)Jan Ob{\l}{\'o}j]{Obloj:2004}
{Ob{\l}{\'o}j, Jan}.
\newblock The {S}korokhod embedding problem and its offspring.
\newblock \emph{Probab. Surv.}, 1:\penalty0 321--390, 2004.
\newblock ISSN 1549-5787.
\newblock \doi{10.1214/154957804100000060}.

\bibitem[Pedersen and Peskir(2001)Jesper~L. Pedersen and Goran
  Peskir]{PedersenPeskir:2001}
{Pedersen, Jesper L. and Peskir, Goran}.
\newblock The {A}z\'ema-{Y}or embedding in non-singular diffusions.
\newblock \emph{Stochastic Process. Appl.}, 96\penalty0 (2):\penalty0 305--312,
  2001.
\newblock ISSN 0304-4149.
\newblock \doi{10.1016/S0304-4149(01)00120-X}.

\bibitem[Perkins(1986)Edwin~A. Perkins]{Perkins:1985}
{Perkins, Edwin A.}
\newblock The {C}ereteli-{D}avis solution to the {$H^1$}-embedding problem and
  an optimal embedding in {B}rownian motion.
\newblock In \emph{Seminar on stochastic processes, 1985 ({G}ainesville,
  {FL})}, volume~12 of \emph{Progr. Probab. Statist.}, pages 172--223.
  Birkh\"{a}user, Boston, MA, 1986.
\newblock \doi{10.1007/978-1-4684-6748-2_12}.

\bibitem[Root(1969)David~H. Root]{Root:1969}
{Root, David H.}
\newblock The existence of certain stopping times on {B}rownian motion.
\newblock \emph{Ann. Math. Statist.}, 40\penalty0 (2):\penalty0 715--718, 1969.
\newblock ISSN 0003-4851.
\newblock \doi{10.1214/aoms/1177697749}.

\bibitem[R\"ost(1971)Hermann R\"ost]{Rost:1971}
{R\"ost, Hermann}.
\newblock The stopping distributions of a {M}arkov {P}rocess.
\newblock \emph{Invent. Math.}, 14\penalty0 (1):\penalty0 1--16, 1971.
\newblock ISSN 0020-9910.
\newblock \doi{10.1007/BF01418740}.

\bibitem[R\"ost(1976)Hermann R\"ost]{Rost:1976}
{R\"ost, Hermann}.
\newblock Skorokhod stopping times of minimal variance.
\newblock In \emph{S\'eminaire de {P}robabilit\'es {X} ({P}remi\`ere partie,
  {U}niv. {S}trasbourg, {S}trasbourg, ann\'ee universitaire 1974/75)}, volume
  511 of \emph{Lecture Notes in Math.}, pages 194--208. Springer, Berlin, 1976.
\newblock ISBN 978-3-540-07681-0.
\newblock \doi{10.1007/BFb0101107}.

\bibitem[Skorohod(1965)Anatoliy~V. Skorohod]{Skorokhod:1965}
{Skorohod, Anatoliy V.}
\newblock \emph{Studies in the theory of random processes}.
\newblock \emph{Adiwes International Series in Mathematics, 7021}.
  Addison-Wesley Publishing Co., Inc., Reading, MA, 1965, pages viii+199.
\newblock ISBN 978-0-201-07021-7.
\newblock Translated from the Russian by Scripta Technica, Inc.

\bibitem[Vallois(1983)Pierre Vallois]{Vallois:1983}
{Vallois, Pierre}.
\newblock Le probl\`eme de {S}korokhod sur {${\bf R}$}: une approche avec le
  temps local.
\newblock In \emph{S\'eminaire de {P}robabilit\'es {XVII} (1981/82)}, volume
  986 of \emph{Lecture Notes in Math.}, pages 227--239. Springer-Verlag, Berlin
  Heidelberg, 1983.
\newblock ISBN 978-3-540-12289-0.
\newblock \doi{10.1007/BFb0068320}.

\bibitem[Wang(2011)Jiajie Wang]{J.Wang:Thesis}
{Wang, Jiajie}.
\newblock \emph{{R}oot's and {R}ost's embeddings: construction , optimality and
  applications to variance options}.
\newblock PhD thesis, University of Bath, 2011.
\newblock \sc url: \tt\url{http://opus.bath.ac.uk/42009/}\normalfont.

\end{thebibliography}
	
\end{document}